\documentclass[fleqn,reqno,11pt,a4paper,final]{amsart}

\usepackage[a4paper,left=20mm,right=20mm,top=20mm,bottom=20mm,marginpar=20mm]{geometry}
\usepackage{amsmath}
\usepackage{amssymb}
\usepackage{amsthm}
\usepackage{amscd}
\usepackage[utf8]{inputenc}
\usepackage{bbm}
\usepackage{xcolor}
\usepackage{color}
\usepackage[english=american]{csquotes}
\usepackage[final]{graphicx}
\usepackage{hyperref}
\usepackage{calc}
\usepackage{mathptmx}
\usepackage{comment}
\usepackage{mathtools}
\usepackage{tikz}
\usepackage{graphicx}
\usepackage{pgfplots}
\pgfplotsset{compat=1.18}
\usetikzlibrary{angles}
\usepackage{enumitem}
\usepackage{stix}

\usepackage[backend=biber,bibencoding=utf8,firstinits,maxnames=4,style=alphabetic,isbn=false, doi=false, eprint= true, url= false,date=year]{biblatex}
\definecolor{luh-dark-blue}{rgb}{0.0, 0.313, 0.608}
\definecolor{luh-light-blue}{rgb}{0.6, 0.725, 0.847}

\graphicspath{{../Pictures/}}

\numberwithin{equation}{section}

\newtheoremstyle{thmlemcorr}{10pt}{10pt}{\itshape}{}{\bfseries}{.}{10pt}{{\thmname{#1}\thmnumber{ #2}\thmnote{ (#3)}}}
\newtheoremstyle{thmlemcorr*}{10pt}{10pt}{\itshape}{}{\bfseries}{.}\newline{{\thmname{#1}\thmnumber{ #2}\thmnote{ (#3)}}}
\newtheoremstyle{remexample}{10pt}{10pt}{}{}{\bfseries}{.}{10pt}{{\thmname{#1}\thmnumber{ #2}\thmnote{ (#3)}}}
\newtheoremstyle{ass}{10pt}{10pt}{}{}{\bfseries}{.}{10pt}{{\thmname{#1}\thmnumber{ A#2}\thmnote{ (#3)}}}

\theoremstyle{thmlemcorr}
\newtheorem{theorem}{Theorem}
\numberwithin{theorem}{section}
\newtheorem{lemma}[theorem]{Lemma}

\newtheorem{proposition}[theorem]{Proposition}

\newtheorem{definition}[theorem]{Definition}
\newtheorem{assumption}[theorem]{Assumption}

\theoremstyle{thmlemcorr*}
\newtheorem*{theorem*}{Theorem}
\newtheorem{lemma*}[theorem]{Lemma}
\newtheorem{corollary*}[theorem]{Corollary}
\newtheorem{proposition*}[theorem]{Proposition}
\newtheorem{problem*}[theorem]{Problem}
\newtheorem{conjecture*}[theorem]{Conjecture}
\newtheorem{definition*}[theorem]{Definition}
\newtheorem{assumption*}[theorem]{Assumption}

\theoremstyle{remexample}
\newtheorem{remark}[theorem]{Remark}
\newtheorem{example}[theorem]{Example}

\theoremstyle{ass}

\newcommand{\Lcal}{\mathcal{L}}

\newcommand{\Zbb}{\mathbb{Z}}

\DeclareMathOperator{\supp}{supp}

\newcommand{\ee}{\mathrm{e}}
\newcommand{\ii}{\mathrm{i}}
\newcommand{\set}[2]{\left\{\, #1 \ \textup{:}\ #2 \,\right\}}

\newcommand{\setb}[2]{\bigl\{\, #1 \ \textup{:}\ #2 \,\bigr\}}
\newcommand{\setB}[2]{\Bigl\{\, #1 \ \textup{:}\ #2 \,\Bigr\}}

\newcommand{\norm}[1]{\|#1\|}

\newcommand{\normb}[1]{\bigl\|#1\bigr\|}

\newcommand{\normBB}[1]{\biggl\|#1\biggr\|}

\newcommand{\abs}[1]{|#1|}
\newcommand{\abslr}[1]{\left|#1\right|}

\newcommand{\absb}[1]{\bigl|#1\bigr|}

\newcommand{\N}{\mathbb{N}}
\newcommand{\R}{\mathbb{R}}
\renewcommand{\S}{\mathbb{S}}

\newcommand{\Z}{\mathbb{Z}}

\newcommand{\eps}{\varepsilon}
\newcommand{\lrangle}[1]{\langle #1 \rangle}

\def\XXint#1#2#3{{\setbox0=\hbox{$#1{#2#3}{\int}$}
\vcenter{\hbox{$#2#3$}}\kern-.5\wd0}}

\renewcommand{\eps}{\varepsilon}
\renewcommand{\phi}{\varphi}

\begin{document}


\title[]{A second-order relaxation system for a degenerate-parabolic quasilinear fourth-order rimming-flow equation with critical spectrum}

\author{Rahul Barthwal}
\address{\textit{Rahul Barthwal:} Institute of Applied Analysis and Numerical Simulation, University of Stuttgart, Pfaffenwaldring~57, 70569 Stuttgart, Germany}
\email{rahul.barthwal@mathematik.uni-stuttgart.de}

\author{Juri Joussen}
\address{\textit{Juri Joussen:} Institute of Analysis, Dynamics and Modeling, University of Stuttgart, Pfaffenwaldring~57, 70569 Stuttgart, Germany}
\email{juri.joussen@iadm.uni-stuttgart.de}

\author{Christina Lienstromberg}
\address{\textit{Christina Lienstromberg:}  Institute of Analysis, Dynamics and Modeling, University of Stuttgart, Pfaffenwaldring~57, 70569 Stuttgart, Germany}
\email{christina.lienstromberg@iadm.uni-stuttgart.de}

\author{Christian Rohde }
\address{\textit{Christian Rohde:} Institute of Applied Analysis and Numerical Simulation, University of Stuttgart, Pfaffenwaldring~57, 70569 Stuttgart, Germany}
\email{christian.rohde@mathematik.uni-stuttgart.de}

\begin{abstract}
For a spectrally critical nonlinear degenerate-parabolic PDE of fourth-order, describing a capillary-driven fluid flow in a rotating cylinder, we propose a second-order relaxation system of parabolic-elliptic type. For positive initial data, we prove local existence and uniqueness of solutions to the relaxation system. Thereby, the lifespan of the approximate solutions turns out to be independent of the relaxation parameter. Moreover, we prove convergence of non-negative weak solutions of the relaxation system to weak solutions of the original equation on compact time intervals. We classify positive steady states of the relaxed system, prove a stability result yielding global strong solutions, and show that the long-time dynamics remain consistent under the relaxation limit.
\end{abstract}
\vspace{4pt}

\maketitle


\bigskip





\section{Introduction}
In this article, we introduce a second-order relaxation system for a fourth-order thin-film rimming-flow equation and analyse this system with respect to the existence and uniqueness of solutions and their convergence to solutions of the rimming-flow equation, as the relaxation parameter tends to infinity. A \emph{rimming flow} is the flow of a fluid covering the inner wall of a rotating horizontal cylinder. 
In our case, we deal with incompressible viscous Newtonian fluids that occupy a thin fluid layer at the impermeable solid inner cylinder wall, and we assume that the fluid flow is driven by surface tension. In order to simplify the problem, we consider only a cross-section of the cylinder and assume that the fluid flow is homogeneous in the direction of the cylinder axis, cf. Figure \ref{fig:1d_rimming-flow} for a sketch of the setting. The resulting equation
\begin{equation}\label{eq:PDE_intro}
    \begin{cases}
        u_t + \big(m(u)(u_{xxx} + u_x)\big)_x = 0,
        &
        t > 0,\ x \in \S^1,
        \\
        u(0,\cdot) = u_0, & x \in \S^1
    \end{cases}
\end{equation}
is a quasilinear degenerate-parabolic partial differential equation of fourth order for the evolution of the height $u$ of the fluid film. More precisely, in \eqref{eq:PDE_intro} the unknown $u=u(t,x) > 0$ describes the positive height of the fluid film at time $t \geq 0$ and angle $x \in \S^1 \simeq\R\ (\mathrm{mod}\ 2\pi)$, measured from the cylinder wall. The non-negative function $m\in C^\infty(\R)$, satisfying $m\geq 0$ and $m(s) = 0$ if and only if $s = 0$, is a the so-called mobility function. Its explicit form is determined by the choice of the boundary condition at the liquid-solid interface. A no-slip condition
entails $m(u)=u^3$, while the Navier slip condition leads to $m(u) = u^3 + \lambda u^2$, where $\lambda >0$ is the slip length.



\begin{center}
\begin{figure}[h]
\begin{tikzpicture}[xscale=0.5, yscale=0.5]
	\draw [thick, fill=luh-dark-blue!50] circle (5);
	\fill [fill=white] (0:4.25) 
	[out=90,in=300] to (30:4) 
	[out=120,in=330]   to (60:4.25) 
	[out=150,in=0]  to (90:4) 
	[out=180,in=30] to (120:4.25) 
	[out=210,in=60]  to (150:4) 
	[out=240,in=90] to (180:4.25) 
	[out=270,in=120] to (210:4) 
	[out=300,in=150] to (240:4.25) 
	[out=330,in=180] to (270:4) 
	[out=0,in=210] to (300:4.25) 
	[out=30,in=240] to (330:4) 
	[out=60,in=270] to (0:4.25);
	\filldraw circle (0.05);
	\draw [thick,->,color=magenta] (5.5,0) arc (0:270:5.5) 
    {}; 
	\draw [thick] (0,0) -- (0:5) {};
    \draw [thick] (0,0) -- (30:5) {};
    \draw [thick,->] (2,0) arc (0:30:2) node[midway,label=right:$x$] {};
	\draw [thick] (210:5) -- (210:4) node[midway,label=right:$\varepsilon$] {}; 
\end{tikzpicture}
\caption{The 1d rimming-flow problem}
\label{fig:1d_rimming-flow}
\end{figure}
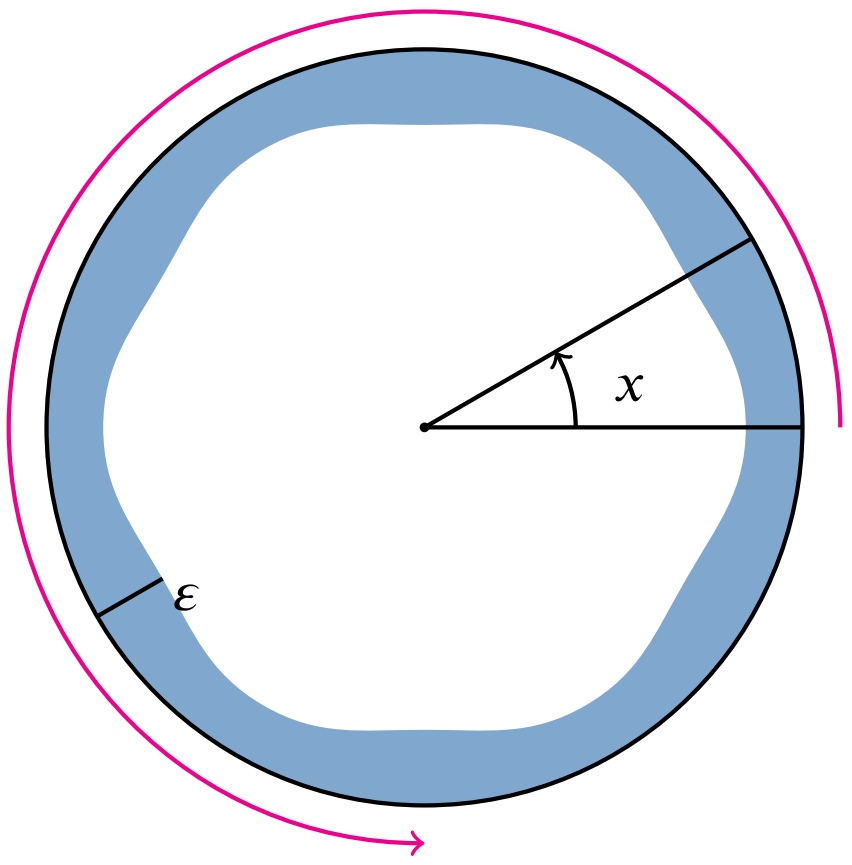 
\end{center}


The underlying equations for incompressible viscous rimming flow are the incompressible Navier--Stokes equations with free surface $u$. This paper, however, concerns the physical situation of 
\begin{itemize}
    \item the fluid layer being rather thin, i.e.
    $$
        \eps =\frac{\text{average film height}}{\text{cylinder radius}} \ll 1;
    $$
    \item the fluid flow being capillary-driven. For the precise scaling assumptions we refer to \cite[Equations (2-8)--(2-11)]{jlv26}.
\end{itemize}
Consequently, \eqref{eq:PDE_intro} is derived from the incompressible Navier--Stokes system with a certain capillary-driven scaling in the formal
asymptotic limit $\eps \to 0$. 

Rimming flows arise in a variety of technical applications such as paper drying or rotational moulding \cite{white1956residual,yih1960instability,Throne_Gianchandani_1980}. In general, the underlying fluid dynamics can exhibit a rich interplay between various competing effects like
rotation, gravity, viscosity, surface tension, and geometric effects. A comprehensive overview, containing additional applications and further physically relevant quantities is given in \cite{Seiden_Thomas_2011}. The mathematical investigation of rimming-flows (and coating flows, where the fluid layer is on the outer wall of the rotating cylinder) starts with \cite{Moffatt_1977,pukhnachev1978motion}. In \cite{Moffatt_1977}, the author studies the gravity-driven equation
\begin{equation}\label{eq:PDE_Moffat}
    u_t + \big(u - g u^3 \cos(x)\big)_x = 0, 
    \quad
    t > 0,\, x \in \S^1,
\end{equation}
while in \cite{pukhnachev1978motion} also surface tension is taken into account, and the corresponding equation is of the form
\begin{equation}\label{eq:PDE_Puk}
    u_t + \big(u + \gamma u^3 (u_{xxx} + u_x) - g u^3 \cos(x)\big)_x = 0, 
    \quad
    t > 0,\, x \in \S^1.
\end{equation}
A detailed derivation and rigorous analysis of this equation is presented in \cite{pukhnachev1978motion,jlv26}. The analogous three-dimensional case is treated in \cite{pukhnachov2005capillary,jllv2026}. A first numerical investigation of the underlying free-surface Navier--Stokes problem is given in \cite{orr1978rimming}. Further literature on thin-film rimming-flow equations similar to \eqref{eq:PDE_Moffat} and \eqref{eq:PDE_Puk}, containing mostly formal analytic and numerical results, has been provided over the last decades. We refer to \cite{ashmore_hosoi_stone_2003,benilov2005does,benilov2008steady,burchard2012convergence,chugunova2009nature,chugunova2010nonnegative,karabash_multi-parameter_2024,lopes2018multiple,thiele2011depinning} and the references therein. 

The rimming-flow problem is closely related to the well-known Taylor--Couette problem, where the fluid flow between two rotating cylinders is studied \cite{couette1890distinction,mallock1896experiments,taylor1923stability}. Of particular interest to this work is the analysis in \cite{pernas_castano_analysis_2020}. The authors derive equation \eqref{eq:PDE_intro} under specific scaling assumptions and rigorously investigate its dynamics. A similar analysis for the non-Newtonian case is provided in \cite{lienstromberg_analysis_2022,LV2024}.

We also note the similarity of equation
\eqref{eq:PDE_intro} to the well-known classical thin-film equation
\begin{equation}\label{eq:thin-film}
    u_t + \gamma \big(u^3 u_{xxx}\big)_x = 0,
    \quad t > 0,\ x \in \Omega \subseteq \R,
\end{equation}
on a flat impermeable substrate, first studied in \cite{BF1990}. In general, however, thin-film rimming-flow equations involve more terms due
to the circular geometry, the rotation and the effects of further involved physical quantities. Indeed, the rotation leads, via a no-slip boundary condition for viscous fluids to a linear transport term $u_x$. In the present paper, we consider the co-moving frame, i.e. the case where fluid is transported exactly at the speed at which the cylinder rotates around its horizontal axis. Thus, using the coordinate transform $x \mapsto x-t$, the term $u_x$ in \eqref{eq:PDE_Moffat} and \eqref{eq:PDE_Puk} does not appear in \eqref{eq:PDE_intro}. Moreover, we assume that surface tension is the driven force, giving rise to the fourth-order term $\big(m(u)(u_{xxx} + u_x)\big)_x$. Other physical effects like inertia or gravity are not taken into account in this paper.

When compared to the classical thin-film equation \eqref{eq:thin-film}, the rimming-flow equation \eqref{eq:PDE_intro} has an additional second-order term $\big(m(u) u_x\big)_x$, a nonlinear backward heat term. On the circle of length $2\pi$, this results in the rimming-flow equation \eqref{eq:PDE_intro} being spectrally critical in the sense that the linearisation around constants $L\phi = -m_0 (\phi_{xxxx} + \phi_{xx})$ with a constant $m_0 > 0$ possesses a non-trivial kernel of first Fourier modes (i.e. sine and cosine functions). On the one hand, this criticality leads to interesting dynamics on the centre manifold, further investigated in Section \ref{sec:long-time-behaviour}. On the other hand, it makes several arguments more involved as the dynamics of the first Fourier modes usually must be treated separately from the remaining (and decaying) Fourier modes. In particular, this means that our results on the rimming-flow equation apply in a natural way to the classical thin-film equation with generally easier proofs.

The main concern of our work is to relax the fourth-order rimming-flow equation \eqref{eq:PDE_intro} by the linear second-order parabolic-elliptic relaxation system 
\begin{equation}\label{eq:relaxed_system}
    \begin{cases}
        v_t^\alpha - \big(m(w^\alpha) \alpha(v^\alpha-w^\alpha)_x\big)_x = 0, 
        \quad t > 0,\ x \in \S^1 &
        \\
        w^\alpha_{xx} + w^\alpha + \alpha(v^\alpha-w^\alpha) = 0, 
        \quad t > 0,\ x \in \S^1 &
        \\
        v^\alpha(0,x) = v_0(x), \quad x \in \S^1.
    \end{cases}
\end{equation}
As the relaxation parameter $\alpha$ tends to infinity, we show that solutions $(v^\alpha,w^\alpha)_{\alpha > 1}$ with $v_0 = u_0$ converge to a solution $u$ of \eqref{eq:PDE_intro}, cf. Theorem \ref{thm:convergence}. Note that we impose an initial condition only on $v^\alpha$, as the initial value $w^\alpha(0,\cdot)$ is then determined by solving $\eqref{eq:relaxed_system}_2$ for $v^\alpha = v_0$. The construction of lower-order approximate systems is often advantageous to gain analytical and numerical insights into nonlinear higher-order PDEs like \eqref{eq:PDE_intro}. 
This strategy is now well established for higher-order PDEs with linear principal parts; see, e.g., \cite{corli2012singular, corli2014parabolic, barthwal2026hyperbolic, dhaouadi2025first}. In particular, Corli and Rohde \cite{corli2012singular} introduced a coupled second-order parabolic–elliptic system to approximate weak solutions of diffusive–dispersive laws. For higher-order PDEs with nonlinear mobilities, Barthwal et al.~\cite{barthwal2026energy} recently proposed a hyperbolic approximation for a broad class of thin-film equations with monotonically increasing disjoining pressure. However, their approximation cannot be applied directly to equations of the form \eqref{eq:PDE_intro}.

Our motivation to introduce the system \eqref{eq:relaxed_system} is twofold. On the one hand, we aim to provide a `nice' approximating system that provides partial linearization and can rather easily be investigated numerically by standard, structure-preserving methods for lower-order systems. Moreover, the relaxation approach paves the way to exploit existing methods for post-processing issues like 
model order reduction or uncertainty quantification. From an analytical perspective, we require this approximation to be well-posed and, moreover, its solutions to converge to solutions of the original fourth-order equation.

On the other hand, the relaxation system \eqref{eq:relaxed_system} can also be viewed as an alternative to classical regularisation schemes for thin-film type equations (cf. \cite{BF1990}). Indeed, we construct weak solutions to the relaxation system \eqref{eq:relaxed_system} that are shown to satisfy an energy equality. In the limit $\alpha \to \infty$, these `regularised' solutions converge to weak solutions of the rimming-flow equation \eqref{eq:PDE_intro} satisfying only an energy inequality.

We point out that our relaxation approach -- while well-suited for numerical analysis -- did not allow us to establish the existence of global non-negative weak solutions by means of an entropy similar to that of Bernis and Friedman \cite{BF1990}. The main obstacle towards such entropy estimates is the nonlocal dependence of $w^\alpha$ on $v^\alpha$ through $\eqref{eq:relaxed_system}_2$, which prevents the cancellation of the mobility factor $m(w^\alpha)$. In that regard, our approach seems to be weaker than the procedure from \cite{BF1990}, which relies on the regularisation of the mobility function.

A similar relaxation model has been proposed in \cite{perthame_relaxation_2021} for Cahn--Hilliard-type equations. Choosing $\psi(n) = -\frac{n^2}{2}$, $n = v^\alpha$, $\sigma = \alpha^{-1}, \lambda = 1$, and $\phi = \alpha (v^\alpha - w^\alpha)$, their approximate system reads (with Neumann-type boundary conditions instead of periodic boundary conditions)
\begin{equation}\label{eq:relaxed_system_perthame}
    \begin{cases}
        v_t^\alpha - \big(m(v^\alpha) \alpha(v^\alpha-w^\alpha)_x\big)_x = 0, 
        \quad t > 0,\ x \in \S^1 &
        \\
        w^\alpha_{xx} + w^\alpha + \alpha(v^\alpha-w^\alpha) = 0, 
        \quad t > 0,\ x \in \S^1.
    \end{cases}
\end{equation}
In this form, the only difference is that the mobility depends on $v^\alpha$ instead of $w^\alpha$ and has two zeros at $0$ and $1$. Using $m(v^\alpha)$ allows them to prove the non-negativity of the mobility by the entropy argument as in \cite{BF1990}. In particular, this implies that $w^\alpha$ is bounded in the interval $[0,1]$. This is crucial in their argument to prove the convergence to solutions of the original Cahn--Hilliard equation. Since in our case we consider a mobility with only one zero, solutions could in general be unbounded, and the convergence proof of \cite{perthame_relaxation_2021} is not applicable. For this reason, we chose $m = m(w^\alpha)$. On the one hand, this has the drawback that the entropy argument of \cite{BF1990} for proving non-negativity of the approximate function is not available. On the other hand, our approach allows us to obtain uniform Hölder bounds for $w^\alpha$ and therewith uniform convergence to the rimming-flow equation. Moreover, given initial data close to a constant, we can prove long-time convergence to a steady state of the relaxation system and the consistency with the long-time behaviour of solutions to the rimming-flow equation: The long-time limit of the relaxed solution converges to the long-time limit of the rimming-flow equation as $\alpha\to\infty$, see Theorem \ref{thm:long-time-consistency}.


\subsection{Notation}

We briefly introduce some notation regarding function spaces and Fourier series used in subsequent sections of this work. For a function $u\in L^2(\S^1)$ we denote its Fourier coefficients by
\begin{equation*}
    \hat u_k \coloneq \fint_{\S^1} u(x)\,\ee^{-\ii k x}\,dx,\quad k\in\Zbb,
\end{equation*}
such that its Fourier series reads
\begin{equation*}
    u = \sum_{k\in\Zbb}\, \hat u_k \ee^{\ii k x}\quad\text{in }L^2(\S^1).
\end{equation*}
Thence, for $s\geq 0$ we define the usual Sobolev spaces by
\begin{equation*}
    H^s(\S^1) \coloneq \setB{u\in L^2(\S^1)}{\sum_{k\in\Zbb}(1 + k^2)^s\abs{\hat u_k}^2 < \infty}
\end{equation*}
with norm
\begin{equation*}
    \norm{u}_{H^s} \coloneq 2\pi \Big(\sum_{k\in\Zbb} (1 + k^2)^s \abs{\hat u_k}^2\Big)^{1/2}.
\end{equation*}
Note that with this definition the $H^0$-norm agrees with the usual $L^2$-norm. Finally, we define $H^{-s}(\S^1)$ as usual to be the dual space of $H^s(\S^1)$, with the operator norm.

It turns out that the rimming-flow equation \eqref{eq:PDE_intro} is spectrally critical with respect to the first Fourier modes. In order to control the first Fourier modes we will therefore sometimes split a function into its first Fourier modes and the remaining Fourier modes. For this reason we define the projections
\begin{equation}\label{eq:def-projection-1}
    P_{\pm 1} u \coloneq \hat u_1 \ee^{\ii x} + \hat u_{-1} \ee^{-\ii x}
\end{equation}
onto the first Fourier modes and
\begin{equation}\label{eq:def-projection-geq2}
    P_{\geq 2} u \coloneq \sum_{\abs{k}\geq 2} \hat u_k\,\ee^{\ii k x}
\end{equation}
onto the Fourier modes greater than or equal to $2$.

Moreover, we define the usual parabolic Hölder spaces
\begin{equation*}
\begin{aligned}
    C^{\rho,\sigma}\big([0,T]\times\S^1\big) \coloneq \big\{u\in & C\big([0,T]\times\S^1\big):
    \abs{u(t_1,x) - u(t_2,x)}\leq C\abs{t_1 - t_2}^\rho \text{ and } \\
    &\abs{u(t,x_1) - u(t,x_2)} \leq C\abs{x_1 - x_2}^\sigma \text{ for all }t,t_1,t_2\in[0,T],\ x,x_1,x_2\in\S^1\big\},
\end{aligned}
\end{equation*}
where, by a slight abuse of notation, we denote by $\abs{x_2-x_1}$ for $x_1,x_2\in\S^1$ the shortest distance between $x_1$ and $x_2$ on $\S^1$.

Finally, by $C > 0$ we denote a generic positive constant which may change from line to line.


\subsection{Main results}
In Section \ref{sec:convergence} we start by defining weak solutions both to the rimming-flow equation \eqref{eq:PDE_intro} and to the relaxation system \eqref{eq:relaxed_system}. Using an energy dissipation mechanism (Lemma \ref{lem:weak_energy_dissipation}), we prove in Lemma \ref{lem:hoelder-continuity} that solutions to the relaxation system \eqref{eq:relaxed_system} are Hölder continuous, uniformly in the relaxation parameter $\alpha\to\infty$. Using the Arzelà--Ascoli Theorem, we extract a uniformly convergent subsequence. This is the key ingredient to prove that weak solutions to the relaxation system \eqref{eq:relaxed_system} converge to weak solutions of the rimming-flow equation \eqref{eq:PDE_intro}, see Theorem \ref{thm:convergence}.

In Section \ref{sec:local-existence} we prove in Theorem \ref{thm:local-well-posedness} local existence and uniqueness of strong solutions to the relaxation system \eqref{eq:relaxed_system} for positive initial data by means of Banach's fixed point theorem. Moreover, for positive initial values $v_0$, the existence interval may be chosen independently of the relaxation parameter $\alpha\to\infty$, thus allowing us to apply the convergence result Theorem \ref{thm:convergence} in Section \ref{sec:convergence}.

Finally, in Section \ref{sec:long-time-behaviour} we investigate the steady states of the relaxation system \eqref{eq:relaxed_system} and its long-time behaviour. In Lemma \ref{lem:steady_states} we classify all positive steady states of the relaxation system. These are given by
\begin{equation*}
\begin{cases}
    v^\alpha(x) = \frac{\alpha-1}{\alpha} a + b \cos(x) + c \sin(x) &
    \\
    w^\alpha(x) = a + b \cos(x) + c \sin(x)
\end{cases}
\end{equation*}
with $b^2 + c^2 < a^2$, corresponding to circular surface profiles possibly shifted from the centre of the cylinder.
In Theorem \ref{thm:stability-relaxed} we then proceed to show that the constant steady states ($b = c = 0$) are stable and that for each initial value close enough to the constant, the solution chooses one steady state (i.e. coefficients $b,c$) to which it converges in the long-time limit. However, the dependence of this long-time limit on the initial data is not yet understood. With Theorem \ref{thm:stability-cited} we cite a similar result for the original rimming-flow equation. We close the paper by proving in Theorem \ref{thm:long-time-consistency} that the long-time dynamics are consistent under the relaxation limit: For a fixed initial value, the long-time limits of the relaxation system converge to the long-time limit of the rimming-flow equation as $\alpha\to\infty$. In future work, this consistency result might allow to further investigate the choice of the long-time limit depending on the initial value by considering the relaxation system \eqref{eq:relaxed_system}.


\section{Convergence of weak solutions of the relaxation system to weak solutions of the limit equation}\label{sec:convergence}

In this section, we prove convergence of weak solutions of the relaxation system
\begin{equation}\label{eq:relaxed_system_conv}
    \begin{cases}
        v_t^\alpha - \big(m(w^\alpha) \alpha(v^\alpha-w^\alpha)_x\big)_x = 0, 
        \quad t > 0,\ x \in \S^1 &
        \\
        w^\alpha_{xx} + w^\alpha + \alpha(v^\alpha-w^\alpha) = 0, 
        \quad t > 0,\ x \in \S^1 &
        \\
        v^\alpha(0,x) = v_0(x), \quad x \in \S^1
    \end{cases}
\end{equation}
with $v_0 = u_0$ to weak solutions of the rimming-flow equation
\begin{equation}\label{eq:rimming-flow_conv}
\begin{cases}
    u_t + \big(m(u) (u_{xxx} + u_x)\big)_x = 0,\quad t>0,\ x\in\S^1 & \\
    u(0,x) = u_0(x), \quad x \in \S^1.
\end{cases}
\end{equation}
Here and in the following we always assume that $m$ is a so-called mobility function with the following properties:
\begin{assumption}\label{as:mobility-assumptions}
    The mobility function $m\in C^\infty(\R^+)\cap C(\R^+_0)$ is non-negative, i.e. $m(s)\geq 0$ for all $s\geq 0$ and satisfies $m(s) = 0$ if and only if $s = 0$.
\end{assumption}
\noindent Prototypical physical examples for the mobility function are given by $m(s) = s^n$, where $n \geq 1$.

In order to prove the convergence, we start by defining weak solutions to the rimming-flow equation \eqref{eq:rimming-flow_conv}. This definition is similar to those defined in \cite{BF1990,chugunova2010nonnegative}, the main difference being that we assume explicitly the non-negativity $u\geq 0$ and the existence of a weak time derivative $u_t\in L^2\big((0,T);H^{-1}(\S^1)\big)$.

\begin{definition}[Weak solutions of rimming-flow equation]\label{def:weak-solution-rimming-flow}
    Let $u_0\in L^2(\S^1)$ with $u_0\geq 0$ almost everywhere. We call a continuous function $u\in C\big([0,T]\times\S^1\big)$ a weak solution to the rimming-flow equation \eqref{eq:rimming-flow_conv} with initial value $u_0$ if
    \begin{enumerate}[label=(\roman*)]
        \item $u$ is non-negative,
        \item $u\in L^2\big((0,T);H^1(\S^1)\big)$,
        \item $u$ is weakly differentiable in time with $u_t\in L^2\big((0,T);H^{-1}(\S^1)\big)$,
        \item $\lim_{t\searrow 0} u(t) = u_0$ in $L^2(\S^1)$,
        \item in the positivity set $P \coloneq \set{(t,x)\in(0,T)\times\S^1}{u(t,x) > 0}$, $u$ is thrice weakly differentiable with respect to $x$, and
        $$m(u)\,u_{xxx} \in L^2(P),$$
    \end{enumerate}
    and for all test functions $\phi\in L^2\big((0,T);H^1(\S^1)\big)$ we have
    \begin{equation}\label{eq:rimming-flow_weak-solution-identity}
        \int_0^T \lrangle{u_t\,,\,\phi}_{H^{-1},H^1}\,dt - \int_P m(u)\, (u_{xxx} + u_x)\,\phi_x\,d(t,x) = 0.
    \end{equation}
\end{definition}

Next, we introduce the notion of a weak solution to the relaxation system \eqref{eq:relaxed_system_conv}. 

\begin{definition}[Weak solution of the relaxation system]\label{def:weak-sol-relaxed}
    Given an initial value $v_0\in L^2(\S^1)$, we call a pair of functions $(v^\alpha, w^\alpha)$ a weak solution to the initial value problem \eqref{eq:relaxed_system_conv} if
    \begin{equation*}
        v^\alpha, w^\alpha \in L^2\big((0,T); H^1(\S^1)\big)\quad\text{with}\quad v_t^\alpha\in L^2\big((0,T);H^{-1}(\S^1)\big),
    \end{equation*}
    such that $w^\alpha \geq 0$ is non-negative almost everywhere in $[0,T]\times\S^1$ and for all test functions $\phi,\psi\in L^2\big((0,T);H^1(\S^1)\big)$ we have
    \begin{equation*}
    \begin{aligned}
        \int_0^T \lrangle{v_t^\alpha\,,\, \phi}_{H^{-1},H^1} \,dt + \int_0^T\int_{\S^1}m(w^\alpha) \alpha (v^\alpha-w^\alpha)_x \phi_x \,dx\,dt &= 0,\\
        - \int_0^T\int_{\S^1} w_x^\alpha \psi_x\,dx\,dt + \int_0^T\int_{\S^1} \big(w^\alpha + \alpha (v^\alpha-w^\alpha)\big)\psi\,dx\,dt &= 0,
    \end{aligned}
    \end{equation*}
    and $v^\alpha$ satisfies the initial condition in the sense that
    \begin{equation*}
        \lim_{t\searrow 0} v^\alpha(t) = v_0\quad\text{in }L^2(\S^1).
    \end{equation*}
\end{definition}

In the following, we collect some properties of weak solutions $(v^\alpha,w^\alpha)$ to the relaxation system \eqref{eq:relaxed_system_conv}. More precisely, we derive further regularity properties, we show that they conserve their integral mean and that they satisfy a (weak) energy-dissipation equality.

We start by deriving some properties of the solution operator for the elliptic equation $\eqref{eq:relaxed_system_conv}_2$. For this lemma, we require the technical assumption $\alpha \geq 2$. In principle, the solution operator is well-defined for $\alpha > 1$ but it is uniformly bounded only for $\alpha$ bounded away from one. Since we consider anyway the limit $\alpha\to \infty$, our assumption does not pose a restriction.

\begin{lemma}[Solution operator for the elliptic equation]\label{lem:solution-operator-2nd}
Let $\alpha\geq 2$. Then, there exists a continuous solution operator $S_\alpha : H^{-1}(\S^1) \to H^1(\S^1)$ such that $ w = S_\alpha v$ is the unique weak solution of
\begin{equation*}
    w_{xx} + w + \alpha(v-w) = 0.
\end{equation*}
Moreover, the operator $S_\alpha$ has the following properties:
\begin{enumerate}[label=(\roman*)]
    \item $S_\alpha$ is continuous as an operator from $H^k(\S^1)$ to $H^{k+2}(\S^1)$ for $k\geq 0$, with
    $$\norm{S_\alpha}_{\Lcal(H^k,H^{k+2})} \leq \alpha;$$
    \item $S_\alpha$ is uniformly bounded in $\alpha$ as an operator from $H^k(\S^1)$ to $H^k(\S^1)$ for all $k\geq 0$:
    $$\norm{S_\alpha}_{\Lcal(H^k;H^k)}\leq \frac{\alpha}{\alpha-1}\leq 2;$$
    \item $S_\alpha$ is self-adjoint, i.e. for all $u,v\in H^{-1}(\S^1)$ we have
    $$\lrangle{S_\alpha u\,,\,v}_{H^1,H^{-1}} = \lrangle{u\,,\,S_\alpha v}_{H^{-1},H^1}.$$
\end{enumerate}
\end{lemma}

\begin{proof}
    For $\alpha > 1$ the operator is given explicitly in Fourier notation by
    \begin{equation}\label{eq:S_alpha_Fourier}
        S_\alpha v = \sum_{n\in\Z} \frac{\alpha}{n^2 + \alpha - 1}\,\hat v_n\,\ee^{\ii n x}.
    \end{equation}

    {\noindent\bfseries (i):}
    Using $\alpha\geq 2$, we have $n^2 + \alpha - 1 \geq n^2 + 1$ and hence
    \begin{equation*}
        \norm{S_\alpha v}_{H^{k+2}}^2 = 2\pi\sum_{n\in\Z} (1 + n^2)^{k+2} \frac{\alpha^2}{(n^2 + \alpha - 1)^2} \abs{\hat v_n}^2 \leq 2\pi\alpha^2 \sum_{n\in\Z} (1 + n^2)^k\abs{\hat v_n}^2 = \alpha^2 \norm{v}_{H^k}^2.
    \end{equation*}

    \medskip
    {\noindent\bfseries (ii):}
    For all $n\in\Z$ and $\alpha\geq 2$ we have
    \begin{equation*}
        \frac{\alpha}{n^2 + \alpha - 1} \leq \frac{\alpha}{\alpha - 1}\leq 2,
    \end{equation*}
    whence it follows that
    \begin{equation*}
        \norm{S_\alpha v}_{H^{k}}^2 = 2\pi\sum_{n\in\Z} (1 + n^2)^{k} \frac{\alpha^2}{(n^2 + \alpha - 1)^2} \abs{\hat v_n}^2 \leq \frac{2\pi\alpha^2}{(\alpha-1)^2} \sum_{n\in\Z} (1 + n^2)^k\abs{\hat v_n}^2 = \frac{\alpha^2}{(\alpha-1)^2} \norm{v}_{H^k}^2.
    \end{equation*}

    \medskip
    {\noindent\bfseries (iii):}
    This follows since $\tfrac{\alpha}{n^2 + \alpha - 1}$ is real-valued and hence we can see in \eqref{eq:S_alpha_Fourier} that $S_\alpha$ has only real eigenvalues.
\end{proof}

Next, we derive some regularity properties of weak solutions $(v^\alpha,w^\alpha)$ to the relaxation system \eqref{eq:relaxed_system_conv}.

\begin{proposition}[Properties of weak solutions to \eqref{eq:relaxed_system_conv}]\label{prop:weak-solution}
    A weak solution $(v^\alpha, w^\alpha)$ as defined in Definition \ref{def:weak-sol-relaxed} is well-defined and has the following properties:
    \begin{enumerate}[label=(\roman*)]
        \item $v^\alpha\in C\big([0,T];L^2(\S^1)\big)$,
        \item $w^\alpha\in L^2\big((0,T); H^3(\S^1)\big)$ with $w_t\in L^2\big((0,T); H^1(\S^1)\big)$,
        \item $w^\alpha\in C\big([0,T];H^2(\S^1)\big)$,
        \item $w_{xx}^\alpha(t) + w^\alpha(t) + \alpha(v^\alpha(t)-w^\alpha(t)) = 0$ for every $0\leq t\leq T$.
    \end{enumerate}
\end{proposition}

\begin{proof}
    {\noindent\bfseries (i):}
    By assumption we have $v^\alpha\in L^2\big((0,T);H^1(\S^1)\big)$ and $v_t^\alpha\in L^2\big((0,T);H^{-1}(\S^1)\big)$. It is well-known that this implies $v^\alpha\in C\big([0,T];L^2(\S^1)\big)$ (see e.g. \cite[Theorem 5.9.3]{evans_partial_2010}).

    \medskip
    {\noindent\bfseries (ii):}
    Using Lemma \ref{lem:solution-operator-2nd} (i), $v^\alpha\in L^2\big((0,T);H^1(\S^1)\big)$ implies $w^\alpha\in L^2\big((0,T);H^3(\S^1)\big)$ and similarly, it follows from $v_t^\alpha\in L^2\big((0,T);H^{-1}(\S^1)\big)$ that $w_t^\alpha\in L^2\big((0,T);H^{1}(\S^1)\big)$.

    \medskip
    {\noindent\bfseries (iii):}
    Given property (ii), this is well-known, cf. the proof of property (i).

    \medskip
    {\noindent\bfseries (iv):}
    This is an immediate consequence of property (iii).
\end{proof}

\begin{remark}
    Proposition \ref{prop:weak-solution} shows in particular that $w^\alpha_0 \coloneq w^\alpha(0)\in H^2(\S^1)$ is well-defined and satisfies
    \begin{equation*}
        w^\alpha_{0,xx} + w^\alpha_0 + \alpha(v_0-w^\alpha_0) = 0\quad \text{in }L^2(\S^1).
    \end{equation*}
\end{remark}

\begin{lemma}[Weak conservation of mean value] \label{lem:conservation_mean}
    A weak solution $(v^\alpha,w^\alpha)$ of the relaxation system \eqref{eq:relaxed_system_conv} in the sense of Definition \ref{def:weak-sol-relaxed} conserves its mean value, i.e. we have
    \begin{equation*}
        \fint_{\S^1} v^\alpha(t,x) dx = \fint_{\S^1} v_0(x) dx \quad\text{and}\quad \fint_{\S^1} w^\alpha(t,x) dx = \fint_{\S^1}w^\alpha_0(x) dx = \frac{\alpha}{\alpha-1}\fint_{\S^1} v_0(x) dx
    \end{equation*}
    for all $0\leq t\leq T$.
\end{lemma}

\begin{proof}
    This follows by using $\phi \equiv \psi \equiv 1$ in Definition \ref{def:weak-sol-relaxed}.
\end{proof}

Note that due to the non-negativity of $w^\alpha \geq 0$ almost everywhere on $[0,T]\times \S^1$ (cf. Definition \ref{def:weak-sol-relaxed}) we even have conservation of mass in the sense that
\begin{equation*}
    \|w^\alpha(t)\|_{L^1(\S^1)}
    =
    \|w^\alpha_0\|_{L^1(\S^1)}
\end{equation*}
for all $0\leq t\leq T$. However, the functions $v^\alpha$ do not necessarily have a sign.
In order to shorten the notation we sometimes use Fourier notation to denote the mean of a solution $v$ by
\begin{equation*}
    \hat v_0 = \fint_{\S^1} v^\alpha(t,x) dx = \fint_{\S^1} v_0(x) dx.
\end{equation*}
Analogously we write $\hat w^\alpha_0$ for the mean of $w^\alpha(t)$, or equivalently of $w^\alpha_0$.

\medskip

In the following lemma we show that weak solutions $(v^\alpha,w^\alpha)$ in the sense of Definition \ref{def:weak-sol-relaxed} satisfy a (weak) energy-dissipation equality for the energy functional
\begin{equation*}
    E_\alpha[v^\alpha,w^\alpha]
    \coloneq 
    \frac{1}{2} \int_{\S^1} |w^\alpha_x|^2 - |w^\alpha|^2 + \alpha |v^\alpha - w^\alpha|^2\, dx + \pi \hat w_0^2.
\end{equation*}
Here, we have added the constant $\pi\hat w_0^2$ to ensure that the energy functional is always non-negative. This can best be seen by using an equivalent Fourier formulation of the energy:
\begin{equation}\label{eq:energy-fourier-representation}
\begin{aligned}
    E_\alpha[v^\alpha,w^\alpha] &= \frac{1}{2} \cdot 2\pi \sum_{k\in\Zbb} (k^2 - 1) \abs{\hat w_k}^2 + \pi \hat w_0^2 + \frac{\alpha}{2}\norm{v^\alpha - w^\alpha}_{L^2}^2\\
    &= \pi \sum_{k\in\Zbb\setminus\{0\}} (k^2 - 1) \abs{\hat w_k}^2 + \frac{\alpha}{2}\norm{v^\alpha - w^\alpha}_{L^2}^2 \geq 0.
\end{aligned}
\end{equation}

\begin{lemma}[Weak energy equality]\label{lem:weak_energy_dissipation}
    Let $(v^\alpha, w^\alpha)$ be a weak solution to \eqref{eq:relaxed_system_conv} as defined in Definition \ref{def:weak-sol-relaxed}. Then, the energy dissipation equality 
    \begin{equation*}
        E_\alpha[v^\alpha,w^\alpha](t) + \alpha^2 \int_0^t\int_{\S^1} m(w^\alpha) (v^\alpha-w^\alpha)_x^2\,dx\,ds = E_\alpha[v_0, w^\alpha_0]
    \end{equation*}
    is satisfied for all $0\leq t\leq T$.
\end{lemma}

\begin{remark}
This lemma is an important consequence of the fact that the relaxation system \eqref{eq:relaxed_system_conv} is indeed regularising (cf. Proposition \ref{prop:weak-solution}). Note that 
special weak solutions are usually constructed as a limit of regularised functions that satisfy an energy \emph{equality} and hence one recovers an energy dissipation \emph{inequality} in the limit. This method goes back to Bernis and Friedman \cite{BF1990}, and has since then been adapted to various situations for related problems. Our approach is quite similar to these constructions in the sense that solutions of the relaxation system are automatically "good enough" to guarantee dissipation.  
\end{remark}

\begin{proof}
    Using e.g. \cite[Theorem 5.9.3]{evans_partial_2010} in combination with Proposition \ref{prop:weak-solution}, we obtain the identities
    \begin{align}
        \label{eq:energy-identity-1}
        &\frac{1}{2}\norm{w_x^\alpha(t)}_{L^2}^2 - \frac{1}{2}\norm{w^\alpha_{0,x}}_{L^2}^2 = \int_0^t\int_{\S^1} w_x^\alpha\,w_{tx}^\alpha\,dx\,ds,\\
        \label{eq:energy-identity-2}
        &\frac{1}{2}\norm{w^\alpha(t)}_{L^2}^2 - \frac{1}{2}\norm{w^\alpha_0}_{L^2}^2 = \int_0^t\int_{\S^1} w^\alpha\,w_t^\alpha\,dx\,ds,\\
        \label{eq:energy-identity-3}
        &\frac{\alpha}{2}\norm{v^\alpha(t) - w^\alpha(t)}_{L^2}^2 - \frac{\alpha}{2}\norm{v_0 - w^\alpha_0}_{L^2}^2 = \int_0^t \lrangle{v_t^\alpha\,,\,\alpha(v^\alpha-w^\alpha)}_{H^{-1},H^1}\,ds - \int_0^t\int_{\S^1} w_t^\alpha \,\alpha(v^\alpha-w^\alpha)\,dx\,ds.
    \end{align}
    Adding \eqref{eq:energy-identity-1} and \eqref{eq:energy-identity-3} and subtracting \eqref{eq:energy-identity-2} yields
    \begin{equation*}
    \begin{aligned}
        &\frac{1}{2}\Big(\norm{w_x^\alpha(t)}_{L^2}^2 - \norm{w^\alpha(t)}_{L^2}^2\Big) + \frac{\alpha}{2}\norm{v^\alpha(t) - w^\alpha(t)}_{L^2}^2 - \int_0^t \lrangle{v_t^\alpha\,,\,\alpha(v^\alpha-w^\alpha)}_{H^{-1},H^1}\,ds\\
        &= \frac{1}{2}\Big(\norm{w^\alpha_{0,x}}_{L^2}^2 - \norm{w^\alpha_0}_{L^2}^2\Big) + \frac{\alpha}{2}\norm{v_0 - w^\alpha_0}_{L^2}^2,
    \end{aligned}
    \end{equation*}
    where we used the definition of a weak solution for $w^\alpha$ with $\psi = w^\alpha_t\,\chi_{[0,t]}$. This is justified by Proposition \ref{prop:weak-solution} (ii).
    Finally, using the definition of a weak solution for $v^\alpha$ with $\phi = \alpha(v^\alpha-w^\alpha)\,\chi_{[0,t]}$ leads to the equation
    \begin{equation*}
        - \int_0^t \lrangle{v_t^\alpha\,,\,\alpha(v^\alpha-w^\alpha)}_{H^{-1},H^1}\,ds = \alpha^2 \int_0^t\int_{\S^1} m(w^\alpha) (v^\alpha-w^\alpha)_x^2\,dx\,ds
    \end{equation*}
    and thus we conclude the result.
\end{proof}

\begin{lemma}[Uniform Hölder-continuity of $w^\alpha$]\label{lem:hoelder-continuity}
    Let $v_0\in H^1(\S^1)$ with $v_0\geq 0$. Let $T>0$ and let $(v^\alpha,w^\alpha)$, $\alpha\geq 2$, be a family of weak solutions of the relaxation system \eqref{eq:relaxed_system_conv} on $[0,T]$, as defined by Definition \ref{def:weak-sol-relaxed} and with initial value $v_0\in H^1(\S^1)$.
    Then, $w^\alpha$ is uniformly Hölder-continuous. More precisely, we have $w\in C^{1/8,1/2}\big([0,T]\times\S^1\big)$ and there exists a constant $C$ independent of $T$ and $\alpha$ such that for all $t, t_1, t_2\in[0,T]$ and all $x,x_1,x_2\in\S^1$ we have
    \begin{equation*}
        \abs{w^\alpha(t,x)}\leq C, \quad \abs{w^\alpha(t_2,x) - w^\alpha(t_1,x)} \leq C\abs{t_2 - t_1}^{1/8}, \quad\text{and}\quad \abs{w^\alpha(t,x_2) - w^\alpha(t,x_2)} \leq C\abs{x_2 - x_1}^{1/2}.
    \end{equation*}
\end{lemma}

Recall that $\abs{x_2-x_1}$ for $x_1,x_2\in\S^1$ denotes the shortest distance between $x_1$ and $x_2$ on $\S^1$.

\begin{remark}
    At first glance, the Hölder-continuity with exponent $\tfrac{1}{8}$ in time seems distinctly non-optimal as one would expect for a second-order equation that the spatial Hölder exponent $\tfrac{1}{2}$ translates to a temporal exponent of $\tfrac{1}{2}\cdot\tfrac{1}{2} = \tfrac{1}{4}$. And indeed, the above estimates can be improved for a fixed $\alpha$. However, the crucial point of Lemma \ref{lem:hoelder-continuity} is the uniformity in $\alpha$. As we will see that solutions to the relaxation system converge to solutions of the fourth-order equation \eqref{eq:rimming-flow_conv}, no better uniform Hölder regularity can be expected.
\end{remark}

\begin{proof}
    We denote by $w^\alpha_0 \coloneq w^\alpha(0) = S_\alpha v_0$ the initial value of $w^\alpha$. We first show that the initial energy $E_\alpha [v_0,w^\alpha_0]$ is uniformly bounded for $\alpha\geq 2$. Indeed, by Lemma \ref{lem:solution-operator-2nd} (ii) we have $\norm{w^\alpha_0}_{H^1} = \norm{S_\alpha v_0}_{H^1} \leq 2\norm{v_0}_{H^1}$. Moreover, we can estimate
    \begin{equation*}
    \begin{aligned}
        \frac{\alpha}{2}\norm{w^\alpha_0 - v_0}_{L^2}^2
        &= \frac{\alpha}{2}\norm{S_\alpha v_0 - v_0}_{L^2}^2\\
        &= 2\pi^2\alpha \sum_{n\in\Z} \Big(\frac{\alpha}{n^2 + \alpha - 1} - 1\Big)^2 \,\abs{\hat v_{0,n}}^2\\
        &= 2\pi^2\alpha \sum_{n\in\Z}\Big(\frac{n^2-1}{n^2 + \alpha - 1}\Big)^2 \,\abs{\hat v_{0,n}}^2.
    \end{aligned}
    \end{equation*}
    Hence, using that $\absb{\frac{n^2-1}{n^2 + \alpha -1}} \leq \frac{n^2 + 1}{\alpha-1}$ and $\absb{\frac{n^2-1}{n^2 + \alpha -1}} \leq 1$ for $\alpha\geq 2$, $n\in\Z$, we obtain
    \begin{equation}\label{eq:initial-w-L2-convergence}
        \frac{\alpha}{2}\norm{w^\alpha_0 - v_0}_{L^2}^2 \leq \frac{2\pi^2\alpha}{\alpha - 1} \sum_{n\in\Z} (1 + n^2)\,\abs{\hat v_{0,n}}^2 \leq \norm{v_0}_{H^1}^2.
    \end{equation}
    Together with the $H^1$-bound on $w^\alpha_0$, this provides us with a bound on the initial energy:
    \begin{equation}\label{eq:initial-energy-bounded}
        E_\alpha[v_0,w^\alpha_0] \leq C,
    \end{equation}
    where the constant $C>0$ depends only on $v_0\in H^1(\S^1)$.
    
    From the Fourier representation of the energy $E_\alpha$ in \eqref{eq:energy-fourier-representation}, we deduce that all but the first Fourier modes of $w^\alpha$ can be controlled by $E_\alpha[v^\alpha,w^\alpha]$. Indeed, using the projection $P_{\geq 2}$ defined in \eqref{eq:def-projection-geq2}, we have
    \begin{equation*}
        \norm{P_{\geq 2}w^\alpha}_{H^1}^2 = 2\pi \sum_{\abs{k}\geq 2} k^2 \abs{\hat w^\alpha_k}^2 \leq  C_1 \pi \sum_{\abs{k}\geq 2} (k^2 - 1) \abs{\hat w_k^\alpha}^2 \leq C_1 E_\alpha [v^\alpha,w^\alpha].
    \end{equation*}
    By Lemma \ref{lem:weak_energy_dissipation} the energy is dissipative:
    \begin{equation} \label{eq:H1_bound}
        \norm{P_{\geq 2}w^\alpha(t)}_{H^1}^2 \leq C_1 E_\alpha[v^\alpha(t),w^\alpha(t)] \leq C_1 E_\alpha[v_0,w^\alpha_0] \leq C_2.
    \end{equation}
    Since $w^\alpha$ is non-negative, all Fourier coefficients are bounded by the mean value of $w^\alpha$. Thanks to the uniform $H^1(\S^1)$-bound on $w^\alpha_0$ and the fact that the mean values are conserved in time (see Lemma \ref{lem:conservation_mean}), it follows that also the mean values of $w^\alpha(t,\cdot)$ are uniformly bounded. Using \eqref{eq:H1_bound}, we obtain $\norm{w^\alpha(t)}_{H^1}\leq C$, and by Morrey's inequality, we immediately obtain
    \begin{equation*}
        \abs{w^\alpha(t,x)}\leq C \quad\text{and}\quad 
        \abs{w^\alpha(t,x_2) - w^\alpha(t,x_1)} \leq C\abs{x_2 - x_1}^{1/2}
    \end{equation*}
    for all $t\in[0,T]$ and all $x,x_1,x_2\in\S^1$.
    
    It remains to prove the uniform Hölder-continuity in time.  
    For this purpose, we adapt the method used in \cite[Proof of Theorem 1]{cornalba2025h}.
    Let $x_0\in\S^1$ and $t_1,t_2\in[0,T]$ be arbitrary. For convenience, we denote in the following $r\coloneq \abs{t_2-t_1}^{1/4}$. Now, we take a function $\phi\in C^\infty(\S^1)$ with the following properties:
    \begin{equation}\label{eq:properties_phi}
        \phi \geq 0, \quad \phi(x) = 0 \text{ for }x\in\S^1\text{ with }\abs{x-x_0}\geq r, \quad \norm{\phi}_{L^1} = 1, \quad \norm{\phi_x}_{L^2} \leq \frac{C}{r^{3/2}}.
    \end{equation}
    With this, we estimate
    \begin{equation*}
    \begin{aligned}
        \abs{w^\alpha(t_2,x_0) - w^\alpha(t_1,x_0)}
        &\leq \abslr{\int_{\S^1}\big(w^\alpha(t_2,x_0) - w^\alpha(t_2,x)\big) \phi(x)\,dx} + \abslr{\int_{\S^1}\big(w^\alpha(t_2,x) - w^\alpha(t_1,x)\big) \phi(x)\,dx}\\
        &+ \abslr{\int_{\S^1}\big(w^\alpha(t_1,x) - w^\alpha(t_1,x_0)\big) \phi(x)\,dx}.
    \end{aligned}
    \end{equation*}
    By the construction of $\phi$ in \eqref{eq:properties_phi}, we may estimate the first and third integral ($i=1,2$) by
    \begin{equation*}
        \abslr{\int_{\S^1}\big(w^\alpha(t_i,x_0) - w^\alpha(t_i,x)\big) \phi(x)\,dx} \leq Cr^{1/2} \norm{\phi}_{L^1} = Cr^{1/2} = C\abs{t_2 - t_1}^{1/8}.
    \end{equation*}
    On the other hand, assuming without loss of generality that $t_1 < t_2$, the second integral can be rewritten as follows:
    \begin{equation*}
    \begin{aligned}
        \abslr{\int_{\S^1}\big(w^\alpha(t_2,x) - w^\alpha(t_1,x)\big) \phi(x)\,dx}
        &=\abslr{\int_{\S^1} S_\alpha\left[v^\alpha(t_2) - v^\alpha(t_1)\right](x)\,\phi(x) \,dx}\\
        &=\abslr{\int_{t_1}^{t_2}\int_{\S^1} S_\alpha v^\alpha_t(t,x)\, \phi(x)\,dx\,dt},
    \end{aligned}
    \end{equation*}
    where we used Lemma \ref{lem:solution-operator-2nd} (i) to swap $S_\alpha$ and the time derivative for the last equality. The self-adjointness from Lemma \ref{lem:solution-operator-2nd} (iii) then implies
    \begin{equation*}
    \begin{aligned}
        &\abslr{\int_{\S^1}\big(w^\alpha(t_2,x) - w^\alpha(t_1,x)\big) \phi(x)\,dx}\\
        &=\abslr{\int_{t_1}^{t_2} \lrangle{v^\alpha_t(t)\,,\,S_\alpha \phi(t)}_{H^{-1},H^1}\,dt}\\
        &=\abslr{\int_{t_1}^{t_2}\int_{\S^1} m(w^\alpha(t,x))\,\alpha\big(v^\alpha(t,x)-w^\alpha(t,x)\big)_x\,S_\alpha\phi_x(x)\,dx\,dt}\\
        &\leq C \norm{m(w^\alpha)}^{1/2}_{C([0,T]\times\S^1)} \normb{\alpha\sqrt{m(w^\alpha)}\, (v^\alpha-w^\alpha)_x} _{L^2([0,T]\times\S^1)} \norm{\phi_x}_{L^2([t_1,t_2]\times\S^1)}.
    \end{aligned}
    \end{equation*}
    We have already established the boundedness of $w^\alpha$, therefore $m(w^\alpha)$ must be bounded, too. The second norm is uniformly bounded by the energy dissipation from Lemma \ref{lem:weak_energy_dissipation}. Finally, by construction of $\phi$ in \eqref{eq:properties_phi}, we have
    \begin{equation*}
        \norm{\phi_x}_{L^2([t_1,t_2]\times\S^1)} = \norm{\phi_x}_{L^2(\S^1)} \sqrt{t_2-t_1} \leq\frac{Cr^2}{r^{3/2}} = Cr^{1/2} = C\abs{t_2-t_1}^{1/8}.
    \end{equation*}
\end{proof}

\begin{remark}
The result on H\"older continuity in time is essentially the same as in \cite[Lemma 2.1]{BF1990} but uses a significantly easier method adapted from \cite[Proof of Theorem 1]{cornalba2025h}. 
\end{remark}

We are now in a position to prove convergence of weak solutions of the relaxation system \eqref{eq:relaxed_system_conv} to weak solutions of the rimming-flow equation \eqref{eq:rimming-flow_conv}.

\begin{theorem}[Convergence to limit equation]\label{thm:convergence}
    Let $u_0\in H^1(\S^1)$ with $u_0\geq 0$. Let $T>0$ and let $(v^\alpha,w^\alpha)$, $\alpha\geq 2$, be a family of weak solutions of the relaxation system \eqref{eq:relaxed_system_conv} on $[0,T]$, as defined in Definition \ref{def:weak-sol-relaxed}, with initial value $v(0) = u_0$.
    Then, there exists a weak solution $u\in C^{1/8,1/2}\big([0,T]\times\S^1\big)$ of the rimming-flow equation \eqref{eq:rimming-flow_conv} as in Definition \ref{def:weak-solution-rimming-flow} with initial value $u(0,x) = u_0(x)$ for all $x\in\S^1$ and a sequence $(\alpha_n)_{n\in\N}$ with $\lim_{n\to\infty}\alpha_n=\infty$ such that
    \begin{equation}\label{eq:convergence-to-weak-solution-1}
        \begin{cases}
            w^{\alpha_n} \longrightarrow u 
            &
            \text{uniformly on }[0,T]\times\S^1,
            \\
            w^{\alpha_n}\, 
            \xrightharpoonup{\hspace{0.2em}*\hspace{0.2em}} 
            u
            &
            \text{weakly-$\ast$ in } L^\infty\big((0,T);H^1(\S^1)\big),
            \\
            v^{\alpha_n} \longrightarrow u
            &
            \text{strongly in } L^\infty\big((0,T);L^2(\S^1)\big),
            \\
            v^{\alpha_n}_t\, \xrightharpoonup{\hspace{0.2em}\phantom{*}\hspace{0.2em}} u_t
            &
            \text{weakly in } L^2\big((0,T);H^{-1}(\S^1)\big).
        \end{cases}
    \end{equation}
    Furthermore, on the (open) positivity set $P \coloneq \set{(t,x)\in(0,T)\times\S^1}{u(t,x) > 0}$ of $u$ we have
    \begin{equation}\label{eq:convergence-to-weak-solution-2}
    \begin{aligned}
        w^{\alpha_n}_{xxx}\, \xrightharpoonup{\hspace{0.2em}\phantom{*}\hspace{0.2em}}
        u_{xxx}\quad \text{weakly in } L^2_{\mathrm{loc}}(P),
    \end{aligned}
    \end{equation}
    where weak convergence in $L^2_{\mathrm{loc}}(P)$ is understood as weak convergence in $L^2(K)$ for every compact subset $K\subset P$.
\end{theorem}

\begin{proof}
    We start by proving the convergence results in \eqref{eq:convergence-to-weak-solution-1} for some function $u$. Then, we proceed to prove the local convergence result in \eqref{eq:convergence-to-weak-solution-2}. Finally, we show that the identified limit function $u$ is indeed a weak solution to \eqref{eq:rimming-flow_conv}.
    
    \noindent \textbf{(i): } By Lemma \ref{lem:hoelder-continuity}, we obtain a uniform bound on $w^\alpha$ in the Hölder space $C^{1/8,1/2}\big([0,T]\times\S^1\big)$. Thence, it follows by the theorem of Arzelà--Ascoli that a subsequence $w^{\alpha_n}$ converges uniformly to some limit function $u\in C^{1/8,1/2}\big([0,T]\times\S^1\big)$. Since by assumption $w^\alpha \geq 0$ for all $\alpha$, $u$ must be non-negative, too. This proves $\eqref{eq:convergence-to-weak-solution-1}_1$.

    Furthermore, by the energy dissipation from Lemma \ref{lem:weak_energy_dissipation} and the boundedness of the first Fourier modes of $w^{\alpha_n}$, it follows that $w^{\alpha_n}$ is uniformly bounded in $L^\infty\big((0,T);H^1(\S^1)\big)$. Hence, by the theorem of Banach--Alaoglu, there exists a weakly-$\ast$ convergent subsequence (which we denote by $w^{\alpha_n}$, too). Since we already know that $w^{\alpha_n}\to u$ uniformly, it is easy to identify the limit as $u$, thus proving $\eqref{eq:convergence-to-weak-solution-1}_2$. In particular, we find that $u\in L^2\big((0,T);H^1(\S^1)\big)$, which is needed for $u$ being a weak solution to the rimming-flow equation \eqref{eq:rimming-flow_conv} in the sense of Definition \ref{def:weak-solution-rimming-flow}.

    We turn to the convergence results for $v^\alpha$: The definition of the energy, cf. \eqref{eq:energy-fourier-representation}, together with the energy dissipation from Lemma \ref{lem:weak_energy_dissipation} allow us to conclude that $v^{\alpha_n}$ converges to $u$ in an appropriate space:
    \begin{equation*}
    \begin{aligned}
        \norm{v^{\alpha_n}(t) - u(t)}_{L^2}^2
        &\leq 2\norm{v^{\alpha_n}(t) - w^{\alpha_n}(t)}_{L^2}^2 + 2\norm{w^{\alpha_n}(t) - u(t)}_{L^2}^2\\
        &\leq \frac{4}{\alpha_n} E_{\alpha_n}[v^{\alpha_n}(t), w^{\alpha_n}(t)] + 2\norm{w^{\alpha_n}(t) - u(t)}_{L^2}^2\\
        &\leq \frac{4C}{\alpha_n} + 2\norm{w^{\alpha_n}(t) - u(t)}_{L^2}^2
    \end{aligned}
    \end{equation*}
    for all $n\in\N$ and all $t\in[0,T]$. Thus, we have proved    \begin{equation}\label{eq:proof-convergence-of-v}
        v^{\alpha_n} \longrightarrow u\text{ in }L^\infty\big((0,T);L^2(\S^1)\big),
    \end{equation}
    which is $\eqref{eq:convergence-to-weak-solution-1}_3$.
    Now, we consider the time derivative $v^\alpha_t\in L^2\big((0,T);H^{-1}(\S^1)\big)$. By Riesz' representation theorem, there exists a function $\phi\in L^2\big((0,T);H^1(\S^1)\big)$ such that for almost every $t\in[0,T]$ we have
    \begin{equation*}
        \norm{v^\alpha_t(t)}_{H^{-1}}^2 = \lrangle{v^\alpha_t(t)\,,\,\phi(t)}_{H^{-1},H^1}\quad\text{with}\quad\norm{\phi(t)}_{H^1} = \norm{v^\alpha_t(t)}_{H^{-1}}.
    \end{equation*}
    Thus, using Definition \ref{def:weak-sol-relaxed} of a weak solution we can estimate
    \begin{equation*}
    \begin{aligned}
        \norm{v^\alpha_t}_{L^2((0,T);H^{-1})}^2
        &= \int_0^T \lrangle{v^\alpha_t(t)\,,\,\phi(t)}_{H^{-1},H^1}\,dt\\
        &= -\int_0^T\int_{\S^1} m(w^\alpha)\alpha(v^\alpha - w^\alpha)_x\,\phi_x\,dx\,dt\\
        &\leq
        \left(\int_0^T\int_{\S^1} m(w^\alpha)\,\abs{\phi_x}^2\,dx\,dt\right)^{1/2}\ 
        \left(\int_0^T\int_{\S^1} m(w^\alpha)\alpha^2\abs{(v^\alpha - w^\alpha)_x}^2\,dx\,dt\right)^{1/2}\\
        &\leq\norm{m(w^\alpha)}_\infty^{1/2} \, \norm{\phi}_{L^2((0,T);H^1)} \, E_\alpha[u_0,w_0^\alpha]^{1/2}\\
        &=\norm{m(w^\alpha)}_\infty^{1/2} \, \norm{v^\alpha_t}_{L^2((0,T);H^{-1})} \, E_\alpha[u_0,w_0^\alpha]^{1/2}.
    \end{aligned}
    \end{equation*}
    The uniform convergence of $w^\alpha$ also provides a uniform bound for $m(w^\alpha)$. Moreover, we have seen in \eqref{eq:initial-energy-bounded} in the proof of Lemma \ref{lem:hoelder-continuity} that the initial energy $E_\alpha[u_0,w^\alpha_0]$ is uniformly bounded. Hence, dividing by $\norm{v^\alpha_t}_{L^2((0,T);H^{-1}(\S^1))}$, we obtain
    \begin{equation*}
        \norm{v^\alpha_t}_{L^2((0,T);H^{-1})} \leq C.
    \end{equation*}
    Since $L^2\big((0,T);H^{-1}(\S^1)\big)$ is a Hilbert space, the theorem of Banach--Alaoglu implies that (after redefining $\alpha_n$ if necessary)
    \begin{equation*}
        v^{\alpha_n}_t\, \xrightharpoonup{\hspace{0.2em}\phantom{*}\hspace{0.2em}} f \quad\text{weakly in }L^2\big((0,T);H^{-1}(\S^1)\big).
    \end{equation*}
    It turns out that this is enough to show that $u$ is weakly differentiable in time with derivative given by $u_t = f \in L^2\big((0,T);H^{-1}(\S^1)\big)$. Indeed, for all test functions $\phi\in C_c^\infty((0,T))$ and $\psi\in C^\infty(\S^1)$ we have by \eqref{eq:proof-convergence-of-v} and Fubini's theorem
    \begin{equation}\label{eq:time-deriv-convergence-to-u-proof}
        \lim_{n\to\infty} \int_0^T\int_{\S^1} v^{\alpha_n}\, \phi' \psi\,dx\,dt = \int_0^T\int_{\S^1} u\,\phi'\,\psi\,dx\,dt = \int_{\S^1} \int_0^T u\,\phi'\,dt\ \psi\,dx.
    \end{equation}
    On the other hand, applying again Fubini's theorem, we may rewrite
    \begin{equation*}
        \int_0^T\int_{\S^1} v^{\alpha_n}\, \phi' \psi\,dx\,dt = \int_{\S^1} \psi \int_0^T v^{\alpha_n}\, \phi'\,dt\,dx.
    \end{equation*}
    The inner time integral can be interpreted as a Bochner integral. Using the definition of a weak time derivative we have
    \begin{equation*}
        \int_0^T v^{\alpha_n}\, \phi'\,dt = -\int_0^T v^{\alpha_n}_t \,\phi\,dt
    \end{equation*}
    and
    \begin{equation*}
        \int_0^T\int_{\S^1} v^{\alpha_n}\, \phi' \psi\,dx\,dt = \Big\langle-\int_0^T v^{\alpha_n}_t \,\phi\,dt\,,\, \psi\Big\rangle_{H^{-1}, H^1} \longrightarrow \Big\langle-\int_0^T f \,\phi\,dt\,,\, \psi\Big\rangle_{H^{-1}, H^1}
    \end{equation*}
    as $n\to\infty$ by the weak convergence of $v^{\alpha_n}_t$. Thence, since $C^\infty(\S^1)\hookrightarrow H^1(\S^1)$ is dense, comparing this to \eqref{eq:time-deriv-convergence-to-u-proof} we find
    \begin{equation*}
        -\int_0^T f \,\phi\,dt = \int_0^T u\,\phi'\,dt \in C(\S^1)\hookrightarrow H^{-1}(\S^1).
    \end{equation*}
    Thus, $u_t = f \in L^2\big((0,T);H^{-1}(\S^1)\big)$, which completes the proof of equation \eqref{eq:convergence-to-weak-solution-1}.
    In particular, we have {$u\in C\big([0,T];L^2(\S^1)\big)$}.
    
    We show that $u$ satisfies the correct initial condition: Thanks to the estimate
    \begin{equation*}
        \frac{\alpha}{2}\norm{w^\alpha_0 - u_0}_{L^2}^2 \leq \norm{u_0}_{H^1}^2,
    \end{equation*}
    see \eqref{eq:initial-w-L2-convergence} in the proof of Lemma \ref{lem:hoelder-continuity}, we obtain
    \begin{equation*}
        \norm{w^\alpha_0 - u_0}_{L^2} \leq \frac{C}{\sqrt{\alpha}}\norm{u_0}_{H^1} \longrightarrow 0.
    \end{equation*}
    Hence, we find
    \begin{equation*}
        u(0) = \lim_{n\to\infty} w^{\alpha_n}(0) = \lim_{n\to\infty} w^{\alpha_n}_0 = u_0\quad\text{in }L^2(\S^1).
    \end{equation*}

    \noindent \textbf{(ii): } We turn to the local convergence result in \eqref{eq:convergence-to-weak-solution-2}. Let $\phi\in C_c^\infty(P)$ be a smooth test function and denote by $K\coloneq \supp( \phi)$ the support of $\phi$. Since $K\subset P$ is compact, $u$ (and therefore also $m(u)$) is bounded away from zero on $K$. Since $w^{\alpha_n}$ converges to $u$ uniformly, there exists an $\eps>0$ and an $n_0\in\N$ such that $m(w^{\alpha_n})\geq\eps>0$ for all $n\geq n_0$.
    Thence, it follows from the energy dissipation in Lemma \ref{lem:weak_energy_dissipation} and the boundedness of the first Fourier modes of $w^{\alpha_n}$ that
    \begin{equation*}
        \norm{w^{\alpha_n}_{xxx}}_{L^2(K)} \leq \norm{w^{\alpha_n}_{xxx} + w^{\alpha_n}_x}_{L^2(K)} + \norm{w^{\alpha_n}_x}_{L^2([0,T]\times\S^1)} = \alpha\,\norm{(v^\alpha-w^\alpha)_x}_{L^2(K)} + \norm{w^{\alpha_n}_x}_{L^2([0,T]\times\S^1)}\leq C_{\eps,T}.
    \end{equation*}
    Thus, (a subsequence of) $w^{\alpha_n}_{xxx}$ converges weakly in $L^2(K)$ to $u_{xxx}\in L^2(K)$:
    \begin{equation*}
        \int_K w^{\alpha_n}_{xxx} \,\phi\, d(t,x) = \int_P w^{\alpha_n}_{xxx} \,\phi\,d(t,x) = -\int_P w^{\alpha_n}\,\phi_{xxx}\,d(t,x) \longrightarrow - \int_P u\,\phi_{xxx}\,d(t,x).
    \end{equation*}
    In particular, $u$ is thrice weakly differentiable with respect to $x$ in $P$. As the test function $\phi$ was arbitrary, this also proves the local weak convergence result $\eqref{eq:convergence-to-weak-solution-2}$.

    \noindent \textbf{(iii): } Finally, we need to show that the limit function $u$ is indeed a weak solution of the rimming-flow equation \eqref{eq:rimming-flow_conv} in the sense of Definition \ref{def:weak-solution-rimming-flow}. The idea of the proof is similar to the convergence proof in \cite[Proposition 3.1]{ansini2004doubly}. We have already verified assumptions (i)--(iv) of Definition \ref{def:weak-solution-rimming-flow}. Therefore, it remains to prove (v), i.e. that
    \begin{equation*}
        m(u)\,u_{xxx} \in L^2(P),
    \end{equation*}
    and the weak solution identity \eqref{eq:rimming-flow_weak-solution-identity}.
    To this end, we use $\eqref{eq:relaxed_system_conv}_2$ to rewrite the nonlinear term as
    \begin{equation*}
        m(w^{\alpha_n})\,\alpha_n(v^{\alpha_n} - w^{\alpha_n})_x = -m(w^{\alpha_n})\,(w^{\alpha_n}_{xxx} + w^{\alpha_n}_x).
    \end{equation*}
    Note that by Proposition \ref{prop:weak-solution} (ii) and (iii) this is a well-defined $L^2\big([0,T]\times\S^1)\big)$ function, although this will in general not be true for the (weak) limit $m(u)(u_{xxx}+u_x)$. We start by observing from the energy dissipation in Lemma \ref{lem:weak_energy_dissipation} that
    \begin{equation}\label{eq:convergence-nonlinear-term-bounded}
        \normb{\sqrt{m(w^{\alpha_n}})\,\alpha_n (v^{\alpha_n} - w^{\alpha_n})_x}_{L^2([0,T]\times\S^1)} \leq C
    \end{equation}
    is uniformly bounded in $\alpha_n$. Again, from the theorem of Banach--Alaoglu it thus follows immediately that there exists a function $g\in L^2\big([0,T]\times\S^1\big)$ such that (after redefining the subsequence)
    \begin{equation}\label{eq:convergence-weak-L2-convergence-to-g}
        \sqrt{m(w^{\alpha_n})}\,(w^{\alpha_n}_{xxx} + w^{\alpha_n}_x)\, \xrightharpoonup{\hspace{0.2em}\phantom{*}\hspace{0.2em}} g\text{ weakly in }L^2\big([0,T]\times\S^1\big).
    \end{equation}
    It remains to identify the limit $g$.
    For an arbitrary function $\phi\in L^2\big([0,T]\times\S^1)\big)$ we write
    \begin{equation}\label{eq:convergence-nonlinear-term}
    \begin{array}{rcl}
        \lefteqn{-\int_0^T\int_{\S^1} m(w^{\alpha_n})\,(w^{\alpha_n}_{xxx} + w^{\alpha_n}_x) \,\phi\,dx\,dt}
        \\[1.7ex]
        &= &\displaystyle -\int_0^T\int_{\S^1} \big(\sqrt{m(w^{\alpha_n})}-\sqrt{m(u)}\big)\sqrt{m(w^{\alpha_n})} \,(w^{\alpha_n}_{xxx} + w^{\alpha_n}_x) \,\phi\,dx\,dt\\[1.7ex]
        &&{}\displaystyle - \int_0^T\int_{\S^1} \sqrt{m(u)}\,\sqrt{m(w^{\alpha_n})} \,(w^{\alpha_n}_{xxx} + w^{\alpha_n}_x) \,\phi\,dx\,dt.
    \end{array}
    \end{equation}
    Since $w^{\alpha_n}\to u$ uniformly, the first integral on the right-hand side converges to zero as $n\to\infty$. For the second integral, we use the weak convergence to $g$ to obtain
    \begin{equation}\label{eq:convergence-g}
        - \int_0^T\int_{\S^1} \sqrt{m(u)}\,\sqrt{m(w^{\alpha_n})} \,(w^{\alpha_n}_{xxx} + w^{\alpha_n}_x) \,\phi\,dx\,dt \longrightarrow - \int_0^T\int_{\S^1} \sqrt{m(u)}\,g\,\phi\,dx\,dt
    \end{equation}
    as $n\to\infty$.
    We take an arbitrary $\eps>0$ and split the domain $[0,T]\times\S^1$ of the integral into
    \begin{equation*}
        P_\eps \coloneq \set{(t,x)\in[0,T]\times \S^1}{m(u(t,x))\geq \eps}
    \end{equation*}
    and its complement $P_\eps^c\coloneq([0,T]\times\S^1)\setminus P_\eps$. Then, on the one hand, using the Cauchy--Schwarz inequality and \eqref{eq:convergence-nonlinear-term-bounded} we obtain
    \begin{equation}\label{eq:convergence-p-eps-c-estimate} \hspace*{-0.2cm}
        \abslr{\int_{P_\eps^c} \sqrt{m(u)}\,\sqrt{m(w^{\alpha_n})} \,(w^{\alpha_n}_{xxx} + w^{\alpha_n}_x) \,\phi\,d(t,x)} 
        \leq C\sqrt{\eps}\,\normb{\sqrt{m(w^{\alpha_n})}\,(w^{\alpha_n}_{xxx} + w^{\alpha_n}_x)}_{L^2([0,T]\times\S^1)} \leq C\sqrt{\eps},
    \end{equation}
    since $\sqrt{m(u)} \leq \sqrt{\eps}$ in $P_\eps^c$. On the other hand, using the assumptions on $m$, $P_\eps\subset P$ is compact and hence $w^{\alpha_n}_{xxx} + w^{\alpha_n}_{x}$ is weakly convergent in $L^2(P_\eps)$. We can identify the limit as $u_{xxx} + u_x$. Since we also have the uniform convergence $\sqrt{m(w^\alpha_n)}\to \sqrt{m(u)}$, we find
    \begin{equation}\label{eq:convergence-p-eps}
        \int_{P_\eps} \sqrt{m(u)}\,\sqrt{m(w^{\alpha_n})} \,(w^{\alpha_n}_{xxx} + w^{\alpha_n}_x) \,\phi\,d(t,x) \longrightarrow \int_{P_\eps} m(u)\,(u_{xxx} + u_x)\,\phi\,d(t,x).
    \end{equation}
    Combining \eqref{eq:convergence-g}, \eqref{eq:convergence-p-eps-c-estimate} and \eqref{eq:convergence-p-eps} we find, as $\eps\to 0$:
    \begin{equation}\label{eq:convergence-definition-g}
        g =
        \begin{cases}
            \sqrt{m(u)}\,(u_{xxx} + u_x) & \text{in }P,\\
            0 & \text{in }\big([0,T]\times\S^1\big) \setminus P.
        \end{cases}
    \end{equation}
    Note that $g\in L^2\big([0,T]\times\S^1\big)$ and since $u$ is bounded and $u_x\in L^2\big([0,T]\times\S^1\big)$, we conclude $m(u)\,u_{xxx} \in L^2(P)$, verifying Definition \ref{def:weak-solution-rimming-flow} (v). Finally, in view of \eqref{eq:convergence-nonlinear-term} and $\eqref{eq:convergence-to-weak-solution-1}_4$, we deduce that $u$ is indeed a weak solution in the sense of Definition \ref{def:weak-solution-rimming-flow}.
\end{proof}

If a weak solution is constructed as a limit of solutions of the relaxation problem \eqref{eq:relaxed_system_conv} by means of Theorem \ref{thm:convergence}, then it satisfies at least an energy dissipation \emph{inequality}. This is in line with the classical results obtained by Bernis and Friedman \cite{BF1990} and is made precise by the following lemma.

\begin{lemma}[Weak energy dissipation inequality]
    Let $u\in C^{1/8,1/2}\big([0,T]\times\S^1\big)$ be a weak solution of the rimming-flow equation \eqref{eq:rimming-flow_conv} in the sense of Definition \ref{def:weak-solution-rimming-flow}. We define the energy functional 
    \begin{equation*}
        E_\infty[u] \coloneq \frac{1}{2} \int_{\S^1} \abs{u_x}^2 - \abs{u}^2 \,dx + \pi\hat u_0^2.
    \end{equation*}
    If $u$ is a limit as constructed in Theorem \ref{thm:convergence}, then it satisfies  energy dissipation inequality
    \begin{equation*}
        E_\infty[u](t) + \int_{P_t} m(u) (u_{xxx} + u_{x})^2\,d(s,x) \leq E_\infty[u_0],
    \end{equation*}
    where $P_t \coloneq P \cap \big((0,t)\times\S^1\big)$ is the positivity set up to time $t$, with $P$ being defined as in Theorem \ref{thm:convergence}.
\end{lemma}

\begin{proof}
    We first show that
    \begin{equation}\label{eq:convergence-initial-energy}
        \lim_{\alpha\to\infty} E_\alpha [u_0,w^\alpha_0] = E_\infty[u_0].
    \end{equation}
    For this purpose, we need to prove that
    \begin{equation}\label{eq:energy-ineq-convergence-results}
        \lim_{\alpha\to\infty} \norm{w^\alpha_0 - u_0}_{H^1} = 0
        \quad\text{and}\quad
        \lim_{\alpha\to\infty} \frac{\alpha}{2} \norm{u_0 - w^\alpha_0}_{L^2}^2 = 0.
    \end{equation}
    Indeed, we have
    \begin{equation*}
        \norm{w^\alpha_0 - u_0}_{H^1}^2 = 4\pi^2 \sum_{n\in\Z} (1+n^2) \Big(\frac{\alpha}{n^2+\alpha-1}-1\Big)^2 \,\abs{u_{0,n}}^2 = 4\pi^2 \sum_{n\in\Z} (1+n^2) \Big(\frac{n^2-1}{n^2+\alpha-1}\Big)^2 \,\abs{u_{0,n}}^2.
    \end{equation*}
    Since $\absb{\frac{n^2-1}{n^2+\alpha-1}}\leq 1$ and $u_0\in H^1(\S^1)$, it follows from dominated convergence (Tannery's theorem) that
    \begin{equation*}
        \lim_{\alpha\to\infty} \norm{w^\alpha_0 - u_0}_{H^1}^2 = 4\pi^2\sum_{n\in\Z} \lim_{\alpha\to\infty}\Big[(1+n^2) \Big(\frac{n^2-1}{n^2+\alpha-1}\Big)^2 \,\abs{u_{0,n}}^2\Big] = 0.
    \end{equation*}
    In order to prove the second result of \eqref{eq:energy-ineq-convergence-results}, observe that
    \begin{equation*}
        \frac{\alpha}{2} \norm{u_0 - w^\alpha_0}_{L^2}^2 = 4\pi^2 \sum_{n\in\Z} \frac{\alpha}{2} \, \Big(\frac{n^2 - 1}{n^2+\alpha-1}\Big)^2\,\abs{u_{0,n}}^2.
    \end{equation*}
    Similar to \eqref{eq:initial-w-L2-convergence}, we estimate $\absb{\frac{n^2-1}{n^2 + \alpha -1}} \leq \frac{n^2 + 1}{\alpha-1}$ and $\absb{\frac{n^2-1}{n^2 + \alpha -1}} \leq 1$ for $\alpha\geq 2$, $n\in\Z$ to obtain
    \begin{equation*}
        \frac{\alpha}{2} \, \Big(\frac{n^2 - 1}{n^2+\alpha-1}\Big)^2\,\abs{u_{0,n}}^2 \leq \frac{\alpha}{2(\alpha-1)}(1+n^2)\leq 1+n^2
    \end{equation*}
    for all $n\in\Z$ and $\alpha\geq 2$. Since $u_0\in H^1(\S^1)$, we may use again dominated convergence and find
    \begin{equation*}
        \lim_{\alpha\to\infty} \frac{\alpha}{2} \norm{u_0 - w^\alpha_0}_{L^2}^2 = 4\pi^2 \sum_{n\in\Z} \lim_{\alpha\to\infty} \left[\frac{\alpha}{2} \, \Big(\frac{n^2 - 1}{n^2+\alpha-1}\Big)^2\,\abs{u_{0,n}}^2\right] = 0.
    \end{equation*}
    Thus, we have proved \eqref{eq:energy-ineq-convergence-results} and, as an immediate consequence, also \eqref{eq:convergence-initial-energy}.
    
    We proceed with the treatment of the dissipation term. Using the weak lower semi-continuity of norms, we deduce from \eqref{eq:convergence-weak-L2-convergence-to-g} with \eqref{eq:convergence-definition-g} that
    \begin{equation}\label{eq:weak-lower-semicontinuity-of-dissipation}
        \int_{P_t} m(u) (u_{xxx} + u_{x})^2\,d(s,x) \leq \liminf_{n\to\infty} \ \alpha_n^2 \int_0^t\int_{\S^1} m(w^{\alpha_n}) (v^{\alpha_n}-w^{\alpha_n})_x^2\,dx\,ds
    \end{equation}
    for all $t \in [0,T]$.
    We know from Lemma \ref{lem:weak_energy_dissipation} that the energy $E_\alpha[v^\alpha,w^\alpha](t)$ is bounded uniformly in $\alpha$ for every $t \in [0,T]$. In particular, there exists a subsequence $w^{\alpha_n}(t)$ which is weakly convergent in $H^1(\S^1)$ for every $t \in [0,T]$. By the uniform convergence of $w^{\alpha_n}$ from Theorem \ref{thm:convergence} we can identify the weak limit as $u(t)$. Hence, we conclude
    \begin{equation}\label{eq:weak-lower-semicontinuity-of-energy}
    \begin{aligned}
        E_\infty[u](t)
        &= \frac{1}{2}\int_{\S^1}\abs{u_x}^2 - \abs{u}^2\,dx + \pi\hat u_0^2\\
        &\leq \liminf_{n\to\infty} \frac{1}{2}\int_{\S^1}\abs{w^{\alpha_n}_x}^2 - \abs{w^{\alpha_n}}^2\,dx + \pi(\hat w^{\alpha_n}_0)^2\\
        &\leq \liminf_{n\to\infty} E_{\alpha_n}[v^{\alpha_n},w^{\alpha_n}](t)
    \end{aligned}
    \end{equation}
    for all $t\in[0,T]$.
    Finally, inserting \eqref{eq:convergence-initial-energy}, \eqref{eq:weak-lower-semicontinuity-of-dissipation} and \eqref{eq:weak-lower-semicontinuity-of-energy} into the energy dissipation \emph{equality} in Lemma \ref{lem:weak_energy_dissipation} for the relaxed problem, we obtain the desired result.
\end{proof}

\section{Local existence and uniqueness of strong solutions for the relaxation system \eqref{eq:relaxed_system_conv}}\label{sec:local-existence}

In order to prove local existence and uniqueness of strong solutions for positive initial data, 
we need 
as an auxiliary result 
the following maximum principles for the elliptic equation for $w^\alpha$
in the relaxation system \eqref{eq:relaxed_system_conv}.
We state them in the following lemma as properties of the solution operator $S_\alpha$ introduced in Lemma \ref{lem:solution-operator-2nd}:

\begin{lemma}[Maximum principles for $S_\alpha$]\label{lem:maximum-principles}
    Let $\alpha \geq 2$. In addition to Lemma \ref{lem:solution-operator-2nd} the solution operator $S_\alpha$, solving the equation
    \begin{equation*}
        -w_{xx} + (\alpha - 1) w = \alpha v
    \end{equation*}
    by $w = S_\alpha v$, has the following properties when perceived as an operator $S_\alpha\colon H^1(\S^1) \to H^3(\S^1)$:
    \begin{enumerate}
        \item[(i)] $S_\alpha$ preserves non-negativity: If $v\in H^{1}(\S^1)$ is non-negative, i.e. $v(x)\geq 0$ for all $x\in\S^1$, then $w = S_\alpha v\geq 0$ is non-negative in $\S^1$, too.
        \item[(ii)] $S_\alpha$ obeys a strong maximum principle: If $v\in H^{1}(\S^1)$ is non-negative and not identically zero, then $w = S_\alpha v > 0$ is strictly positive in $\S^1$.
        \item[(iii)] $S_\alpha$ preserves positivity uniformly in $\alpha$: If $v\in H^{1}(\S^1)$ is bounded away from zero by a constant $c_0 > 0$, i.e. $v(x)\geq c_0 > 0$ for all $x \in \S^1$, then $w(x) = S_\alpha v(x) \geq \tfrac{\alpha}{\alpha-1} c_0> c_0>0$ for all $x\in\S^1$.
    \end{enumerate}
\end{lemma}

\begin{proof}
    Properties (i) and (iii) follow immediately from the observation that $w = S_\alpha v \in H^3(\S^1)\hookrightarrow C^2(\S^1)$ and hence, if the minimum of $w$ is attained in $x_0\in\S^1$, we must have 
    \begin{equation*}
        w(x) 
        \geq
        w(x_0) = \frac{\alpha}{\alpha-1} v(x_0) + \frac{1}{\alpha-1}w_{xx}(x_0) \geq \frac{\alpha}{\alpha-1} v(x_0)
    \end{equation*}
    for all $x\in\S^1$ since $w_{xx}(x_0)\geq 0$. To obtain (iii) recall that $\frac{\alpha}{\alpha-1} > 1$ for $\alpha \geq 2$.
    For property (ii) we refer the reader e.g. to \cite[Theorem 6.4.4]{evans_partial_2010}.
\end{proof}

We now turn to local existence and uniqueness of solutions to the relaxation initial-value problem \eqref{eq:relaxed_system_conv}. To this end, let $v_0\in H^2(\S^1)$, $v_0>0$ be a strictly positive initial value, i.e. we assume that the fluid layer on the inner wall of the cylinder is coating the wall completely and there are no dry spots. By Lemma \ref{lem:maximum-principles} (iii), the corresponding initial value $w_0^\alpha \coloneq S_\alpha v_0$ for $w$ is positive, too. This positivity is preserved at least for a short time interval. We prove this rigorously in the following theorem:

\begin{theorem}[Local well-posedness of the relaxation system]\label{thm:local-well-posedness}
    Let $v_0 \in H^2(\S^1)$ with $v_0(x)>0$ for all $x\in\S^1$ be a strictly positive initial value. Then, there exists a time $T_\ast>0$, only depending on $\min_{x\in\S^1} v_0(x)$ and $\norm{v_0}_{H^2}$, such that for all $\alpha\geq 2$, there exists a unique solution $(v^\alpha,w^\alpha)$ with
    \begin{equation}\label{eq:lwp-solution-regularity}
    \begin{cases}
        v^\alpha \in C^\rho\big([0,T_\ast);H^2(\S^1)\big)\cap C^{1+\rho}\big((0,T_\ast);L^2(\S^1)\big),
        \\
        w^\alpha \in C^\rho\big([0,T_\ast);H^4(\S^1)\big)\cap C^{1+\rho}\big((0,T_\ast);H^2(\S^1)\big)
    \end{cases}
    \end{equation}
    of the relaxed initial-value problem \eqref{eq:relaxed_system_conv}. The component $w^\alpha$ is strictly positive in $[0,T_\ast)\times\S^1$.
    
    Moreover, for a fixed $\alpha\geq 2$ the solution can be extended either to a global solution ($T_\alpha = \infty$) or to a time $T_\alpha\in(T_\ast,\infty)$ at which
    \begin{equation}\label{eq:lwp-film-rupture}
        \lim_{t\nearrow T_\alpha} \min_{x\in\S^1} w^\alpha(t,x) = 0,
    \end{equation}
    i.e. the minimum of $w^\alpha(t,\cdot)$ converges to zero.
\end{theorem}

Note that the relaxation system \eqref{eq:relaxed_system_conv} degenerates if $w^\alpha$ becomes zero in finite time. By Lemma \ref{lem:maximum-principles} (ii), the positivity of $v_0$ implies that the initial value $w^\alpha_0 = S_\alpha v_0$ of $w^\alpha$ is strictly positive, too. Theorem \ref{thm:local-well-posedness} shows in particular that the degeneracy is the only reason that can prevent the solution from being global: There is no finite-time blow-up without the minimum of $w^\alpha$ becoming zero. On the other hand, the $v^\alpha$-component of the solution can become negative in finite time, even while $w^\alpha$ is still positive (and thus a unique solution exists). An example of such a situation is constructed in Example \ref{ex:negativity-of-v} below.

\begin{proof}
    We start by proving the existence and uniqueness for some positive time $T>0$ (to be chosen later) that may depend on $\alpha$ by means of Banach's fixed point theorem. To this end, let
    \begin{equation*}
        \eps\coloneq \tfrac{1}{2}\inf_{\alpha\geq 2}\min_{x\in\S^1} m(w_0^\alpha(x)).
    \end{equation*}
    In general the time of existence will depend on $\eps$. Due to Lemma \ref{lem:maximum-principles} (iii) and the positivity of $v_0$, $w_0^\alpha$ is bounded away from zero and so is $m(w_0^\alpha)$. Therefore, we have $\eps>0$. Let $M_1, M_2$ be two constants that will be chosen later depending on $\alpha$, $\eps$ and $v_0$. Then, we define the complete metric space
    \begin{equation*}
    \begin{aligned}
        X \coloneq \setb{w \in C^\rho\big([0,T]; H^3(\S^1)\big)}{&\norm{w}_{C([0,T]; H^3)}\leq M_1,\ \norm{w}_{C^\rho([0,T]; H^3)}\leq M_2,\\ & w(0,x) = w_0^\alpha(x),\ m(w(t,x))\geq\eps \text{ for all }t\in[0,T],\ x\in\S^1}
    \end{aligned}
    \end{equation*}
    with the metric induced by the Hölder norm on $C^\rho\big([0,T]; H^3(\S^1)\big)$. In particular, for every $\alpha\geq 2$ we have $w_0^\alpha\in X$ as a function constant in time if $M_1, M_2 \geq \norm{w_0}_{H^3}$.
    
    For an arbitrary $w\in X$, we solve the first equation of the relaxation system $\eqref{eq:relaxed_system_conv}_1$: 
    \begin{equation}\label{eq:lwp-ivp}
    \begin{cases}
        v_t - \alpha \big(m(w) v_x\big)_x = f,
        & 0 < t < T,\ x\in\S^1,
        \\
        v(0,x) = v_0(x),
        & x\in\S^1,
    \end{cases}
    \end{equation}
    where $f\coloneq -\alpha (m(w) w_x)_x \in C^\rho\big([0,T];H^1(\S^1)\big)$. For a fixed $w\in X$ this is a simple inhomogeneous linear parabolic equation, whence by standard theory (see e.g. \cite[Theorem II.1.2.1]{amann_linear_1995}) there exists a unique solution
    \begin{equation}\label{eq:lwp-v-regularity}
        v\in C^\rho\big([0,T]; H^2(\S^1)\big) \cap C^{1+\rho}\big((0,T); L^2(\S^1)\big),
    \end{equation}
    the time $T$ of existence being dependent on $\eps$.
    We denote this solution by $\Phi(w) \coloneq v$.
    Thanks to e.g. \cite[Theorem II.5.2.1]{amann_linear_1995}, $v$ is uniformly bounded in $H^1(\S^1)$ by
    \begin{equation}\label{eq:lwp-selfmap-bounded-v}
        \norm{v}_{C([0,T];H^1)} \leq C\,\ee^{\nu T}\,\left(\norm{v_0}_{H^1} + \sqrt{T}\,\norm{f}_{C([0,T];L^2)}\right) \leq C\,\ee^{\nu T}\,\left(\norm{v_0}_{H^1} + \sqrt{T}\,C(\alpha, M_1)\right).
    \end{equation}
    Moreover, we control its Hölder seminorm by \cite[Theorem II.5.3.1]{amann_linear_1995}
    \begin{equation}\label{eq:lwp-selfmap-hoelder-v}
        \norm{v(t) - v(s)}_{H^1} \leq C\,\ee^{\nu T}\, \abs{t-s}^\rho\,\left(\norm{v_0}_{H^{1+2\rho}} + \norm{f}_{C([0,T];L^2)}\right).
    \end{equation}
    Then, we obtain from \eqref{eq:lwp-selfmap-bounded-v} and \eqref{eq:lwp-selfmap-hoelder-v} the uniform Hölder estimate
    \begin{equation}\label{eq:lwp-uniform-hoelder-v}
        \norm{v}_{C^\rho([0,T];H^1)} \leq C\,\ee^{\nu T}\,\left(\norm{v_0}_{H^{1+2\rho}} + \norm{f}_{C([0,T];L^2)}\right) \leq C\,\ee^{\nu T}\,\left(\norm{v_0}_{H^{1+2\rho}} + C(\alpha, M_1)\right).
    \end{equation}
    
    We want to show that, for appropriately chosen constants $T, M_1, M_2$, the operator $S_\alpha\circ\,\Phi:X\to X$ is a well-defined self-map. Here, we understand that $S_\alpha$ is applied to a time-dependent function in a pointwise manner. To this end, we proceed to solve the second equation
    \begin{equation*}
        -w_{xx} + (\alpha - 1) w = \alpha v
    \end{equation*}
    for a given arbitrary function $v$ satisfying \eqref{eq:lwp-selfmap-bounded-v} and \eqref{eq:lwp-uniform-hoelder-v}. We do so by applying the bounded linear solution operator $S_\alpha$ from Lemma \ref{lem:solution-operator-2nd}. Therewith, we obtain the estimates
    \begin{equation}\label{eq:lwp-selfmap-bounded}
        \norm{w}_{C([0,T];H^3)} \leq C(\alpha)\,\ee^{\nu T}\,\left(\norm{v_0}_{H^1} + \sqrt{T}\,C(\alpha, M_1)\right)
    \end{equation}
    from \eqref{eq:lwp-selfmap-bounded-v}, and
        \begin{equation}\label{eq:lwp-uniform-hoelder}
        \norm{w}_{C^\rho([0,T];H^3)} \leq C(\alpha)\,\ee^{\nu T}\,\left(\norm{v_0}_{H^{1+2\rho}} + C(\alpha, M_1)\right)
    \end{equation}
    from \eqref{eq:lwp-uniform-hoelder-v}. If $T>0$ is small enough, we may choose a constant $M_1$ large enough, depending on $\alpha$ and $\norm{v_0}_{H^1}$, such that
    \begin{equation*}
        \norm{w}_{C([0,T];H^3)} \leq M_1.
    \end{equation*}
    Then, we choose $M_2$ large enough, depending on $\alpha$, $\norm{v_0}_{H^{1+2\rho}}$ and $M_1$, such that
    \begin{equation}\label{eq:lwp-uniform-hoelder-w}
        \norm{w}_{C^\rho([0,T];H^3)} \leq M_2.
    \end{equation}
    Moreover, from \eqref{eq:lwp-uniform-hoelder-w}, Sobolev embeddings and the assumptions on $m$, we may deduce that $m(w^\alpha(t,x))\geq \eps$ for all $t\in[0,T]$ and $x\in\S^1$, provided $T > 0$ is small enough.
    Hence, choosing $T, M_1, M_2$ as indicated above yields $S_\alpha(\Phi(w))\in X$ for all $w\in X$, and we have shown that $S_\alpha\circ\,\Phi:X\to X$ is a well-defined self-map.

    In order to be in a position to apply Banach's fixed point theorem, it remains to prove that $S_\alpha\circ\,\Phi:X\to X$ is a contraction, after possibly choosing $T>0$ even smaller. For this purpose, let $w_1,w_2\in X$ be arbitrary and for $i=1,2$ denote by $f_i = -\alpha\big(m(w_i)w_{i,x}\big)_x$ the corresponding right-hand sides of \eqref{eq:lwp-ivp} and by $v_i \coloneq\Phi(w_i)$ the corresponding solutions. Then, $\tilde v\coloneq v_1 - v_2$ solves
    \begin{equation*}
        \tilde v_t - \alpha \big(m(w_1) \tilde v_x\big)_x = f_1 - f_2 + \alpha\big(\big(m(w_1)-m(w_2)\big) v_{2,x}\big)_x, \quad 0 < t < T,\ x\in\S^1,
    \end{equation*}
    with $\tilde v(0) = 0$. Applying estimate \eqref{eq:lwp-uniform-hoelder-v} to this inhomogeneous parabolic equation for $\tilde v$ yields
    \begin{equation}\label{eq:lwp-hoelder-difference}
        \norm{\tilde v}_{C^\rho([0,T];H^1)} \leq C\,\ee^{\nu T}\,\left(\norm{f_1 - f_2}_{C([0,T];L^2)} + \normb{\alpha\big(\big(m(w_1)-m(w_2)\big) v_{2,x}\big)_x}_{C([0,T];L^2)}\right).
    \end{equation}
    We estimate the right-hand side of \eqref{eq:lwp-hoelder-difference} in terms of $\|w_1 - w_2\|_{C([0,T];H^1)}$. Due to the regularity of $m$, we have
    \begin{equation}\label{eq:lwp-lipschitz-mobility}
        \norm{m(w_1) - m(w_2)}_{C([0,T];H^1)} \leq C\norm{w_1 - w_2}_{C([0,T];H^1)},
    \end{equation}
    and \cite[Theorem II.5.2.1]{amann_linear_1995} implies    \begin{equation}\label{eq:lwp-boundedness-v2}
        \norm{v_2}_{C([0,T];H^2)} \leq C \left(\norm{v_0}_{H^2} + C(M_1)\right).
    \end{equation}
    Moreover, using \eqref{eq:lwp-lipschitz-mobility}
    \begin{equation}\label{eq:lwp_boundedness_f1-f2}
        \|f_1 - f_2\|_{C([0,T];L^2)}
        \leq 
        C \|w_1 - w_2\|_{C([0,T];H^2)}.
    \end{equation} 
    Hence, plugging \eqref{eq:lwp-lipschitz-mobility}, \eqref{eq:lwp-boundedness-v2} and \eqref{eq:lwp_boundedness_f1-f2} into \eqref{eq:lwp-hoelder-difference},
    we find
    \begin{equation}\label{eq:lwp-lipschitz-v}
        \norm{\tilde v}_{C^\rho([0,T];H^1)} \leq C(\alpha,M_1)\norm{w_1 - w_2}_{C([0,T];H^2)}.
    \end{equation}
    On the other hand, since $w_1(0) = w_0 = w_2(0)$, we have for every $0\leq t\leq T$
    \begin{equation}\label{eq:lwp-difference-w-small}
        \norm{w_1(t) - w_2(t)}_{H^2} = \norm{w_1(t) - w_2(t) - (w_1(0) - w_2(0))}_{H^2} \leq T^\rho\norm{w_1 - w_2}_{C^\rho([0,T];H^2)}.
    \end{equation}
    Finally, we use Lemma \ref{lem:solution-operator-2nd} (i), \eqref{eq:lwp-lipschitz-v} and \eqref{eq:lwp-difference-w-small} to conclude
    \begin{equation*}
    \begin{aligned}
        \norm{S_\alpha(\Phi(w_1)) - S_\alpha(\Phi(w_2))}_{C^\rho([0,T];H^3)}
        &\leq C(\alpha) \norm{\Phi(w_1) - \Phi(w_2)}_{C^\rho([0,T];H^1)}\\
        &\leq C(\alpha,M_1) T^\rho \norm{w_1 - w_2}_{C^\rho([0,T]; H^2)}.
    \end{aligned}
    \end{equation*}
    Thus, we can ensure that $S_\alpha\circ\,\Phi$ is indeed a contraction on $X$ by choosing $T$ small enough. Banach's fixed point theorem therefore asserts the existence and uniqueness of a fixed point $w^\alpha\in X$ and thus of a local solution $(v^\alpha, w^\alpha) = (\Phi(w^\alpha),w^\alpha)$ of the relaxed initial-value problem \eqref{eq:relaxed_system_conv}. The regularity \eqref{eq:lwp-solution-regularity} is a direct consequence of Lemma \ref{lem:solution-operator-2nd} (i) and \eqref{eq:lwp-v-regularity}.
    
    Note that a-priori the time of existence $T>0$ may depend on $\alpha$, $\min_{x\in\S^1} v_0(x)$ (via $\eps$) and $\norm{v_0}_{H^2}$. We show in the following that solutions can always be extended at least to some time $T_\ast > 0$ independent of $\alpha$. To this end, let $(v^\alpha,w^\alpha)$, $w^\alpha > 0$, be a maximally extended solution of the initial-value problem \eqref{eq:relaxed_system_conv} and denote by $T_\alpha > 0$ its maximal time of existence. In terms of $T_\alpha$, we have so far shown that $T_\alpha > 0$ for every $\alpha\geq 2$. We will hence prove that
    \begin{equation}\label{eq:lwp-uniform-existence-time}
        T_\ast \coloneq\inf_{\alpha\geq 2} T_\alpha > 0.
    \end{equation}
    In the case $T_\alpha = \infty$ for a certain $\alpha\geq 2$ no constraint is imposed on the choice of $T_\ast$. Therefore, we assume without loss of generality that $T_\alpha < \infty$.

    In order to prove \eqref{eq:lwp-uniform-existence-time}, we need to first show \eqref{eq:lwp-film-rupture}. We start by considering a fixed $\alpha\geq 2$ and assuming for contradiction that $w^\alpha$ is bounded away from zero in $[0,T_\alpha)\times\S^1$, i.e. that in particular $m(w^\alpha(t,x))\geq c_0 > 0$ for all $0\leq t < T_\alpha$ and $x\in\S^1$. We show that this positivity assumption effectively prevents a norm blow-up. To this end, we first point out that for instance by a bootstrapping argument, solutions $(v^\alpha,w^\alpha)$ are actually smooth in their existence interval (cf. also \cite[Chapter 8.3]{lunardi_1995}). Therefore, we may estimate for all $0<t<T_\alpha$:
    \begin{equation}\label{eq:lwp-no-norm-blowup}
    \begin{aligned}
        \frac{d}{dt}\frac{1}{2}\norm{v^\alpha_{xx}(t)}_{L^2}^2
        = & -\int_{\S^1} v^\alpha_{xxx} \big(m(w^\alpha) \alpha (v^\alpha - w^\alpha)_x\big)_{xx} dx\\
        \leq & -c_0 \,\alpha\int_{\S^1}\abs{v^\alpha_{xxx}}^2 \,dx + \alpha \int_{\S^1} v^\alpha_{xxx}\big(m(w^\alpha) w^\alpha_x\big)_{xx}\,dx\\
        &- \alpha \int_{\S^1} v^\alpha_{xxx} \big(m(w^\alpha)\big)_{xx}\, v^\alpha_x\,dx - 2\alpha\int_{\S^1} \big(m(w^\alpha)\big)_x\,v^\alpha_{xx}\,v^\alpha_{xxx}\,dx.
    \end{aligned}
    \end{equation}
    Thanks to the energy dissipation from Lemma \ref{lem:weak_energy_dissipation}, we know that
    \begin{equation*}
        \frac{1}{2\alpha}\norm{w_{xx}^\alpha + w^\alpha}_{L^2}^2 = \frac{\alpha}{2}\norm{v^\alpha-w^\alpha}_{L^2}^2 \leq E_\alpha[v_0,w^\alpha_0],
    \end{equation*}
    and since $w^\alpha$ is non-negative in $[0,T_\alpha)$ by assumption, this yields an $H^2$-bound on $w^\alpha$:
    \begin{equation}\label{eq:lwp-H2-bound}
        \norm{w^\alpha(t)}_{H^2}\leq  C, \quad 0\leq t<T_\alpha,
    \end{equation}
    where $C=C({\norm{v_0}_{H^2}})$.\\
    
    Since $v_0\in H^2(\S^1)$, this bound is indeed even uniform in $\alpha$, although this fact is of no particular significance for this step of the proof. Instead we point out that \eqref{eq:lwp-H2-bound} provides in particular bounds on $\norm{w^\alpha_x}_{L^\infty}$ and $\norm{m(w^\alpha)}_{H^2}$. Hence, we can estimate the right-hand side of \eqref{eq:lwp-no-norm-blowup} by
    \begin{equation*}
    \begin{aligned}
        \frac{d}{dt}\frac{1}{2}\norm{v^\alpha_{xx}(t)}_{L^2}^2
        &\leq -c_0\,\alpha\,\norm{v^\alpha_{xxx}}_{L^2}^2 + C\,\alpha\, \norm{v^\alpha_{xxx}}_{L^2}\big(\norm{w^\alpha_{xxx}}_{L^2} + \norm{v^\alpha_{xx}}_{L^2} + \norm{v^\alpha_{x}}_{L^2} + 1\big)\\
        & = -c_0\,\alpha\,\norm{v^\alpha_{xxx}}_{L^2}^2 + C\,\alpha\, \norm{v^\alpha_{xxx}}_{L^2}\big(\norm{S_\alpha v^\alpha_{xxx}}_{L^2} + \norm{v^\alpha_{xx}}_{L^2} + 1\big)\\
        &\leq -c_0\,\alpha\,\norm{v^\alpha_{xxx}}_{L^2}^2 + C_\alpha\,\norm{v^\alpha_{xxx}}_{L^2}\big(\norm{v^\alpha_x}_{L^2} + \norm{v^\alpha_{xx}}_{L^2} + 1\big) \\
        &\leq -c_0\,\alpha\,\norm{v^\alpha_{xxx}}_{L^2}^2 + C_\alpha\,\norm{v^\alpha_{xxx}}_{L^2}\big(\norm{v^\alpha_{xx}}_{L^2} + 1\big).
    \end{aligned}
    \end{equation*}
    From here, we first use Young's inequality to absorb $\norm{v^\alpha_{xxx}}_{L^2}$ into the dissipative term and then employ Gronwall's inequality to obtain
    \begin{equation*}
        \norm{v^\alpha_{xx}(t)}_{L^2}^2 \leq \big(1 + \norm{v_{0,xx}}_{L^2}^2\big)\,\ee^{C(\alpha,c_0)\,t} \leq C(\norm{v_0}_{H^2},\alpha,c_0,T_\alpha)<\infty,\quad 0\leq t < T_\alpha.
    \end{equation*}
    However, if this were true, then restarting the relaxed initial-value problem \eqref{eq:relaxed_system_conv} from $\tilde v_0\coloneq v^\alpha(t)$, the Banach argument from above would yield an extension by a time $\tilde T > 0$ which (although depending on $\alpha$) will not depend on $t<T_\alpha$. This poses a contradiction to the assumed maximality of $T_\alpha$. Thus, we can conclude that $m(w^\alpha)$ cannot be bounded away from zero, and indeed that
    \begin{equation*}
        \liminf_{t\nearrow T_\alpha}\min_{x\in\S^1} w^\alpha(t,x) = 0.
    \end{equation*}
    Finally, the Hölder continuity of $w^\alpha$ from Lemma \ref{lem:hoelder-continuity} implies that the $\liminf$ is indeed a limit and hence we have proved \eqref{eq:lwp-film-rupture}. In fact, we recall that the Hölder continuity is uniform in $\alpha$. Therefore, since $w^\alpha_0\geq c_0>0$ for some constant $c_0>0$ by Lemma \ref{lem:maximum-principles} (iii), there must exist a time $T_\ast>0$, independent of $\alpha$, such that for all $0<t<\min\{T_\alpha,T_\ast\}$
    \begin{equation*}
        w^\alpha(t,x)\geq \frac{c_0}{2}>0.
    \end{equation*}
    As we have just seen, this already implies $T_\alpha > T_\ast$. This concludes the proof.
\end{proof}

\begin{remark}
    The strong maximum principle from Lemma \ref{lem:maximum-principles} (ii) for the second equation $\eqref{eq:relaxed_system_conv}_2$ of the relaxed problem shows that $w^\alpha_0>0$ in $\S^1$, whenever $v_0\geq 0$, $v_0\not\equiv 0$, is only non-negative but not identically zero. The proof of Theorem \ref{thm:local-well-posedness} then shows that there exists a local strong solution even in this case. However, in contrast to the statement of the theorem, the time of existence might, in this case, in general not be uniform in $\alpha$.
\end{remark}

Theorem \ref{thm:local-well-posedness} shows that the existence of global solutions is equivalent to the $w^\alpha$-component remaining bounded away from zero. For mobility functions of the form $m(s) = s^n$ with $n>3/2$, first numerical experiments seem to indicate that $w^\alpha$ does remain positive and global solutions exist. If that were the case, we could actually use Theorem \ref{thm:convergence} to prove global existence of non-negative weak solutions for the rimming-flow equation \eqref{eq:rimming-flow_conv}. Compared to the usual approach going back to Bernis and Friedman in \cite{BF1990}, this would not rely on a regularisation of the mobility function. However, we have not been able to prove that $w^\alpha$ remains bounded away from zero and thus the global existence of solutions to the relaxation system remains an open problem. Towards this problem, one might have the idea to prove the \emph{non-negativity} of $v^\alpha$ first, and then use the strong maximum principle in Lemma \ref{lem:maximum-principles} (ii) to obtain \emph{positivity} of $w^\alpha$. The following example shows that this strategy cannot be successful by constructing a solution $(v^\alpha,w^\alpha)$, where $v^\alpha$ becomes negative, while $w^\alpha$ remains positive at least for a short time interval:

\begin{example}\label{ex:negativity-of-v}
    Fix $\alpha>1$. While $w^\alpha$ is always positive until its maximum time of existence, we can in general not expect $v$ to remain non-negative, even for strictly positive initial values: Take e.g. $v_0\in C^\infty(\S^1)$ that is non-negative ($v_0\geq 0$) but not identically zero ($v_0\not\equiv 0$) such that $v_0$ vanishes in a neighbourhood of some $x_0\in\S^1$. For $\eps>0$ let then $v_{0,\eps}\coloneq v_0 + \eps$ be a strictly positive initial value.

    By the strong maximum principle, $w^\alpha_0 \coloneq S_\alpha v_0 > 0$ in $\S^1$, hence $v_0$ is an admissible initial value in the sense of the local well-posedness result. Consequently, there exists $T_\alpha > 0$ such that the solution $(v^\alpha,w^\alpha)$ exists and is smooth in $[0,T_\alpha)$. Similarly, let $T_{\alpha,\eps}>0$ denote the maximal time of existence for initial value $v_{0,\eps}$. By the Hölder continuity of $w^\alpha$, we find
    \begin{equation*}
        T_{\alpha,\ast} \coloneq \inf_{\eps> 0} T_{\alpha,\eps} > 0.
    \end{equation*}
    We compute
    \begin{equation*}
    \begin{aligned}
        v_t(0,x_0)
        &= \left.\big(m(w^\alpha_0)\,\alpha (v_0 - w^\alpha_0)_x\big)_x \right\vert_{x = x_0}\\
        &= -\alpha\, m\big(w^\alpha_0(x_0)\big) \,w^\alpha_{0,xx}(x_0) - \alpha\, m'\big(w^\alpha_0(x_0)\big)\,\big(w^\alpha_{0,x}(x_0)\big)^2\\
        &=-\alpha(\alpha-1)\,m\big(w^\alpha_0(x_0)\big)\,w^\alpha_0(x_0) - \alpha\, m'\big(w^\alpha_0(x_0)\big)\,\big(w^\alpha_{0,x}(x_0)\big)^2,
    \end{aligned}
    \end{equation*}
    where we have used $w^\alpha_{0,xx}(x_0) = (\alpha - 1) w^\alpha_0(x_0)$. If the mobility $m$ is increasing, then we obtain $v_t(0,x_0) < 0$, and hence there exists a $t_0\in(0,T_\alpha)\cap(0,T_{\alpha,\ast})$ such that $v^\alpha(t_0,x_0) < 0$. By standard parabolic theory, $v^\alpha_\eps$ converges to $v^\alpha$ uniformly in $[0,t_0]\times\S^1$. Hence, we conclude $v^\alpha_\eps(t_0,x_0) < 0$ for $\eps > 0$ small enough.
\end{example}


\section{Steady states and long-time dynamics}\label{sec:long-time-behaviour}

Thanks to Theorem \ref{thm:local-well-posedness} on the local existence and uniqueness of solutions to the relaxation system \eqref{eq:relaxed_system_conv}, we know that (strong) solutions only cease to exist in finite time when the minimum of $w^\alpha(t,\cdot)$ reaches zero, i.e. if
\begin{equation*}
    \lim_{t\nearrow T_\alpha} \min_{x\in\S^1} w^\alpha(t,x) = 0.
\end{equation*}
Moreover, it is well-known that in general non-negative weak solutions cannot be expected to be unique, see e.g. \cite{BBdP1995,cuesta_self-similar_2018,king_moving_2001} and others. A specific situation admitting global positive strong solutions of the rimming-flow equation \eqref{eq:rimming-flow_conv} is a stability result proved in \cite{pernas_castano_analysis_2020}. The authors show that positive constant steady states of \eqref{eq:rimming-flow_conv} are stable but not asymptotically stable:

\begin{theorem}[{\cite[Theorem 5.4]{pernas_castano_analysis_2020}}]\label{thm:stability-cited}
    Let $\hat{u}_0 > 0$. Then, there exists an $\eps_0>0$ such that for all $0 < \eps\leq \eps_0$ and all initial values $u_0\in H^4(\S^1)$ satisfying
    \begin{equation*}
        \fint_{\S^1} u_0(x)\,dx = \hat{u}_0 \quad\text{and}\quad \norm{u_0 - \hat{u}_0}_{H^4} < \eps,
    \end{equation*}
    the rimming-flow equation \eqref{eq:rimming-flow_conv} admits a unique global positive smooth solution $u\in C^\infty(\R^+\times\S^1)$. Moreover, there exists a long-time limit $u_\ast = \hat{u}_0 + b \cos(x) + c\sin(x)$ with constants $b,c\in\R$ such that
    \begin{equation*}
        \norm{u(t) - u_\ast}_{H^4} \leq C\eps \ee^{-\rho t},
    \end{equation*}
    for constants $\rho = \rho(\hat{u}_0) > 0$ and $C = C(\hat{u}_0) > 0$.
\end{theorem}

We point out that \cite[Theorem 5.4]{pernas_castano_analysis_2020} only states the stability result by itself. The existence of the long-time limit (and thus of the constants $b,c\in\R$) is only contained in the proof of the theorem.

\begin{figure}
    \centering
    \begin{tikzpicture}[scale=1.2]
        \coordinate (M) at (0.0, 0.0);
        \coordinate (Mprime) at (-0.8, 0.0);
        \coordinate (P) at ({-0.8+2*cos(50)}, {2*sin(50)});
        \coordinate (dummy) at (2.0, 0.0);
        \coordinate (Q) at ({(3-1-0.8*0.3021)*0.3021}, {(3-1-0.8*0.3021)*0.9533});

        \draw[black, thick] (M) circle (3);

        \draw[luh-light-blue] (Mprime) -- (P) node[left, pos=0.6] {$1-\hat u_0$};
        \draw[luh-light-blue] (Mprime) -- (M);
        \draw[luh-dark-blue] (M) -- (P) node[right, pos=0.5] {$r$};
        \draw[luh-dark-blue] (M) -- (dummy) node[below, pos=0.55] {$1-\hat u_0$};
        \draw[luh-dark-blue] (P) -- (Q);
        \pic[draw, luh-dark-blue, angle eccentricity=1.0, angle radius=0.7cm] {angle = dummy--M--P} node[anchor=south west, luh-dark-blue, xshift=1.5, yshift=-1.0] {$x$};
        \pic[draw, luh-light-blue, angle eccentricity=1.0, angle radius=0.8cm] {angle = M--Mprime--P} node[anchor=south east, luh-light-blue, xshift=-4.5, yshift=-1.0] {$x'$};
        
        \fill[luh-dark-blue] (M) circle (0.05) node[below] {$M$};
        \fill[luh-light-blue] (Mprime) circle (0.05) node[below] {$M'$};
        \fill[luh-light-blue] (P) circle (0.05) node[left, yshift=-1.5] {$P$};
        \fill[magenta] (Q) circle (0.05);

        \draw[domain=0:360, samples=60, smooth, variable=\theta, luh-dark-blue]
            plot ({(3-1)*cos(\theta)}, {(3-1)*sin(\theta)});
        \draw[domain=0:264, samples=60, smooth, variable=\theta, magenta]
            plot ({(3-1-0.8*cos(\theta))*cos(\theta)}, {(3-1-0.8*cos(\theta))*sin(\theta)})
            node[below] {$u_\ast$}
        [domain=264:360, samples=60, smooth, variable=\theta, magenta]
            plot ({(3-1-0.8*cos(\theta))*cos(\theta)}, {(3-1-0.8*cos(\theta))*sin(\theta)});
        \draw[domain=0:300, samples=60, smooth, variable=\theta, luh-light-blue]
            plot ({-0.8+2*cos(\theta)}, {2*sin(\theta)})
            node[above] {$C'$}
        [domain=300:360, samples=60, smooth, variable=\theta, luh-light-blue]
            plot ({-0.8+2*cos(\theta)}, {2*sin(\theta)});
    \end{tikzpicture}
    \caption{\cite[Figure 4]{jllv2026} Cylinder cross-section (black) with limit profile $u_\ast$ (red) parametrised in $x$. The shifted circle $C'$ around $M'$ (light blue) must be reparametrised in terms of $x$ in order to compare it to the profile $u_\ast$. We assume the distance between $M$ and $M'$ to be small but for illustrational purposes it is exaggerated in this figure.}
    \label{fig:shiftedcircles}
\end{figure}
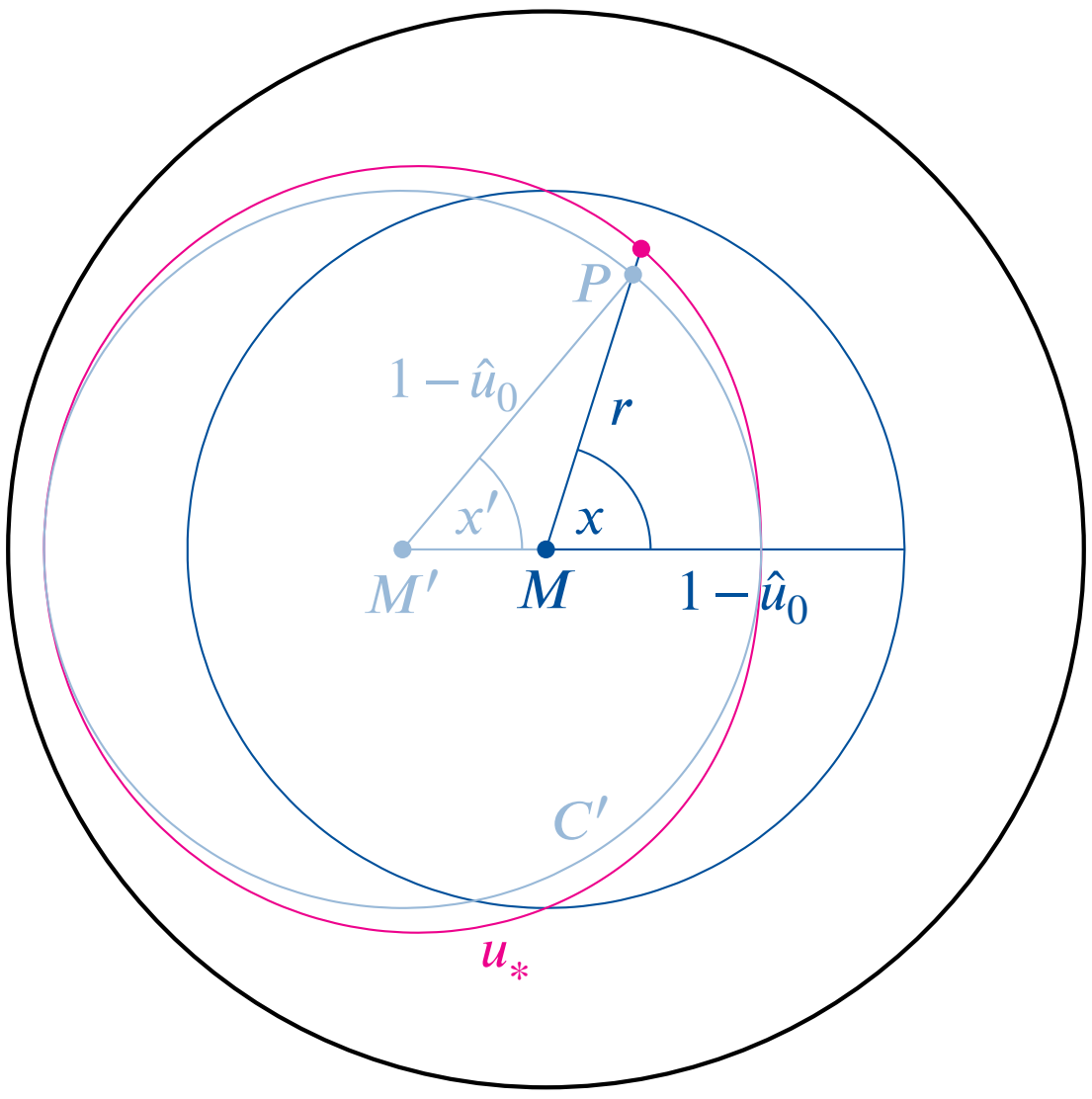

Geometrically, if $b,c\in\R$ are small, steady states of the form $\hat{u}_0 + b \cos(x) + c\sin(x)$ correspond to circular surface profiles whose centre is shifted from the centre of the cylinder by the vector $(-b,-c)$. This can be seen using the law of cosines \cite[Remark 5.6]{jllv2026}. This is illustrated in Figure \ref{fig:shiftedcircles} which is taken with adapted notation from \cite[Figure 4]{jllv2026}.

It is the goal of the present section to prove a result similar to Theorem \ref{thm:stability-cited} for the relaxation system \eqref{eq:relaxed_system_conv}. To this end, we will start by classifying all positive steady states of the relaxation system \eqref{eq:relaxed_system_conv} and then proceed to prove the stability result. 

An obvious question about Theorem \ref{thm:stability-cited} is how the constants $b,c\in\R$ are determined, depending on the initial data. Unfortunately, the theorem is proved by basically showing that the time derivative $u_t$ converges to zero, thus providing no insight into the mechanism choosing $b$ and $c$. If we want to have any hope to approach this question by means of our relaxation system, it is therefore imperative that we show that the long-time limit of the relaxation system converges to the long-time limit of the rimming-flow equation as $\alpha\to\infty$. This is not at all clear from the convergence result in Theorem \ref{thm:convergence}, as its proof relies heavily on the theorem of Arzelà--Ascoli and is hence valid only on any compact (though arbitrarily large) time interval. We will nevertheless prove this consistency in the final part of this section.

\bigskip
\noindent\textbf{\textsc{Characterisation of steady states.}}
In the following we characterise positive steady states of the relaxation system \eqref{eq:relaxed_system_conv} as solutions corresponding to circular surface profiles that may be shifted from the centre of the cylinder.

\begin{lemma}\label{lem:steady_states}
A pair of functions $v^\alpha, w^\alpha \in H^2(\S^1)$ with $w^\alpha(x) > 0,\ x \in \S^1$, is a steady state of the relaxation system \eqref{eq:relaxed_system_conv} if and only if
\begin{equation} \label{eq:steady_states}
\begin{cases}
    v^\alpha(x) = \frac{\alpha-1}{\alpha} a + b \cos(x) + c \sin(x), &
    \\
    w^\alpha(x) = a + b \cos(x) + c \sin(x)
\end{cases}
\end{equation}
for all $x \in \S^1$ and $a, b, c \in \R$ such that $a^2 > b^2 + c^2$.
\end{lemma}

\begin{proof}
\textbf{$\Longleftarrow$:} Let $v, w$ be of the form \eqref{eq:steady_states}. Then $\alpha(v^\alpha-w^\alpha) = -a$ is constant and hence it follows by $\eqref{eq:relaxed_system_conv}_1$ that $v_t^\alpha = 0$ and $w_{xx}^\alpha + w^\alpha = a = -\alpha (v^\alpha-w^\alpha)$.

\medskip
\noindent\textbf{$\Longrightarrow$:}
Let $v^\alpha, w^\alpha \in H^2(\S^1)$ with $m(w^\alpha(x)) > 0,\ x \in \S^1$, be a steady state of \eqref{eq:relaxed_system_conv}. From Lemma \ref{lem:weak_energy_dissipation} we know that
\begin{equation*}
    0
    =
    \frac{d}{dt} E[v^\alpha,w^\alpha]
    =
    -\alpha^2 \int_{\S^1} m(w^\alpha) (v^\alpha-w^\alpha)_x^2\, dx,
\end{equation*}
since $v^\alpha$ and $w^\alpha$ do not depend on time. Using that $m(w^\alpha) > 0$ we conclude that $\alpha (v^\alpha-w^\alpha)_x = 0$. By $\eqref{eq:relaxed_system_conv}_2$ it then follows that $w_{xxx}^{\alpha} + w_x^\alpha = 0$. This implies that
\begin{equation*}
    w^\alpha(x) = a + b \cos(x) + c \sin(x), \quad x \in \S^1,
\end{equation*}
with $a, b, c \in \R$. Again by $\eqref{eq:relaxed_system_conv}_2$ we obtain that $v$ is given by
\begin{equation*}
    v^\alpha(x)
    =
    \frac{\alpha-1}{\alpha} a + b \cos(x) + c \sin(x), 
    \quad x \in \S^1,
\end{equation*}
where $a, b, c \in \R$ have to be chosen in a way such that $w^\alpha > 0$.
\end{proof}

\medskip
\noindent\textbf{\textsc{Global existence of solutions for initial data close to constants.}} In order to obtain the global existence result for initial data close to constants, we need additional higher-order energy estimates. To this end, we introduce for $k \in \N$ the functionals
\begin{equation}\label{eq:higher-order-energy}
    E_\alpha^k[v^\alpha,w^\alpha]
    \coloneq
    \frac{1}{2} \int_{\S^1} |\partial_x^{k+1}w^{\alpha}|^2 - |\partial_x^k w^{\alpha}|^2 + \alpha |\partial_x^k(v^{\alpha}-w^{\alpha})|^2\, dx.
\end{equation}
Similarly to the Fourier representation \eqref{eq:energy-fourier-representation} of the original energy $E_\alpha$, we have
\begin{equation}\label{eq:higher-order-energy-fourier-representation}
    E^k_\alpha[v^\alpha,w^\alpha] = \pi \sum_{\abs{n}\geq 2} n^{2k}(n^2 - 1) \abs{\hat w^\alpha_n}^2 + \frac{\alpha}{2}\norm{\partial_x^k(v^\alpha - w^\alpha)}_{L^2}^2 \geq 0,
\end{equation}
and as such, we can control $\norm{P_{\geq 2} w^\alpha}_{H^k}$ by means of the higher-order energy. We make use of this fact by proving exponential decay of all but the first Fourier modes of solutions. This will then allow us to establish the stability of positive constant steady states:

\begin{theorem}[Stability of constants]\label{thm:stability-relaxed}
The constant steady states obtained in Lemma \ref{lem:steady_states} with $b=c=0$ are stable. More precisely, let $k \geq 3$ and $\hat{v}_0 \in \R$ with $m_0 \coloneq m(\hat{v}_0) > 0$. Then, there exists an $\eps_0 > 0$ such that for all $0<\eps\leq\eps_0$, there exist $\delta = \delta(k,m_0,\eps) > 0$ and $\alpha_0 = \alpha_0(k,m_0,\eps) > 1$ such that for all $\alpha > \alpha_0$ and all initial values $v_0 \in H^{k+1}(\S^1)$ satisfying
\begin{equation*}
    \fint_{\S^1} v_0(x)\, dx = \hat{v}_0
    \quad \text{and} \quad
    \|v_0 - \hat{v}_0\|_{H^{k+1}}
    < \delta,
\end{equation*}
the relaxation system \eqref{eq:relaxed_system_conv} admits a unique global positive smooth solution $v^\alpha, w^\alpha\in C^\infty(\R^+\times \S^1)$. Moreover, there exists a long-time limit
\begin{equation*}
    \begin{cases}
    v^\alpha_\ast(x) = \hat{v}_0 + b_\alpha\cos(x) + c_\alpha \sin(x),\\
    w^\alpha_\ast(x) = \tfrac{\alpha}{\alpha-1}\hat{v}_0 + b_\alpha\cos(x) + c_\alpha \sin(x),
    \end{cases}
\end{equation*}
with constants $b_\alpha, c_\alpha \in \R$ such that
\begin{equation}\label{eq:first-fourier-convergence}
    \|v^\alpha(t) - v_\ast^\alpha\|_{H^k} + \|w^\alpha(t) - w_\ast^\alpha\|_{H^k} \leq C\,\eps\,\ee^{-\rho t},
    \quad t > 0,
\end{equation}
for constants $\rho = \rho(k,m_0) > 0$ and $C=C(k,m_0)>0$.
\end{theorem}

In the proof we need to separately treat all Fourier modes $\geq 2$ and the Fourier modes $0,\pm 1$. This can be seen by considering the linearisation
of the relaxed system around constants $v^\alpha = \hat v_0>0$ with $m_0 \coloneq m(\hat v_0) > 0$. It can be written as
\begin{equation*}
    v^\alpha_t = Lv^\alpha \coloneq -m_0 \big(S_\alpha v^\alpha_{xxxx} + S_\alpha v^\alpha_{xx}\big)
\end{equation*}
with spectrum
\begin{equation*}
    \sigma(L) = \setB{\frac{\alpha(n^2 - n^4)}{n^2 + \alpha - 1}}{n\in\Z},
\end{equation*}
consisting of a critical part $\sigma_0(L) = \{0\}$ corresponding to the non-trivial kernel $\ker L = \mathrm{span}\,\{1,\sin(x),\cos(x)\}$ and an infinite stable part $\sigma_-(L)$, corresponding to the Fourier modes $\geq 2$.

We obtain exponential decay for the Fourier modes $\geq 2$ by showing exponential decay of the higher-order energy $E_\alpha^k$, as defined by \eqref{eq:higher-order-energy}.

In a second step, we project the evolution equation onto the first Fourier modes to obtain
\begin{equation*}
    \frac{d}{dt}\hat v_{\pm 1}(t)\,e^{\pm \ii x} = -P_{\pm 1}\big[\big(m(w^\alpha)(w^\alpha_{xxx}+w^\alpha_x\big)_x\big].
\end{equation*}
Then, we use that the term $(w^\alpha_{xxx} + w^\alpha_x)$ does not contain first Fourier modes and thus decays exponentially. Since the time derivative $\frac{d}{dt} v_{\pm 1}$ decays exponentially, there exists a long-time limit in $\ker L$.

\begin{proof}
According to the above idea, we define the maximal time of closeness
\begin{equation*}
    T_{\eps} 
    \coloneq
    \sup\left\{t > 0\,:\,\|v^\alpha(s) - \hat{v}_0\|_{H^k} \leq \eps,\ s \leq t\right\},
\end{equation*}
until that the solutions is close enough to a constant for the linearisation to be valid.
Note that $T_\eps > 0$ if we choose $\delta \leq \eps$.
We show by contradiction that $T_{\eps} = \infty$ for appropriate initial values, that is, for $\delta > 0$ small enough. To this end, we assume that $T_\eps<\infty$. Then, by the construction of $T_\eps$, we must have
\begin{equation}
    \norm{v^\alpha(T_\eps)-\hat{v}_0}_{H^k} = \eps.
\end{equation}
We derive a differential inequality for the higher-order energy to be able to apply Gronwall's lemma yielding exponential decay of all but the first Fourier modes of $w^\alpha$.
Hence, for $0 < t < T_{\eps}$, we differentiate the higher-order energy $E_\alpha^k[v^\alpha,w^\alpha]$ with respect to time and use $\eqref{eq:relaxed_system_conv}_2$:
\begin{equation*}
\begin{split}
    \frac{d}{dt} E_\alpha^k[v^\alpha,w^\alpha](t)
    &=
    \int_{S^1} \partial_x^{k+1}w^\alpha\,  \partial_x^{k+1}w^\alpha_t - \partial_x^{k}w^\alpha\, \partial_x^{k} w_t
    +
    \alpha \partial_x^{k}(v^\alpha-w^\alpha)\,  \partial_x^{k} (v^\alpha-w^\alpha)_t\, dx
    \\
    &=
    -\int_{S^1} \big(\partial_x^{k+2}w^\alpha + \partial_x^{k}w^\alpha + \alpha \partial_x^{k}(v^\alpha-w^\alpha)\big) \partial_x^{k} w_t\, dx
    +
    \alpha
    \int_{S^1}
    \partial_x^{k}(v^\alpha-w^\alpha)\,  \partial_x^{k} v_t^\alpha\, dx
    \\
    &=
    -\alpha \int_{S^1}
    \partial_x^{k+1}(v^\alpha-w^\alpha)\,  \partial_x^{k-1} v_t^\alpha\, dx.
\end{split}
\end{equation*}
Using $\eqref{eq:relaxed_system_conv}_1$, we find that
\begin{equation*}
\begin{split}
    \frac{d}{dt} & E_\alpha^k[v^\alpha,w^\alpha](t)
    =
    -\alpha^2\int_{\S^1} \partial_x^{k+1}(v^\alpha-w^\alpha)\,  \partial_x^{k}\big(m(w^\alpha) (v^\alpha-w^\alpha)_x\big)\, dx
    \\
    &=
    -\alpha^2\int_{\S^1} \partial_x^{k+1}(v^\alpha-w^\alpha)\,  \partial_x^{k}\big(m_0 (v^\alpha-w^\alpha)_x\big)\, dx
    -
    \alpha^2\int_{\S^1} \partial_x^{k+1}(v^\alpha-w^\alpha)\,  \partial_x^{k}\big((m(w^\alpha)-m_0) (v^\alpha-w^\alpha)_x\big)\, dx
    \\
    &=
    -\alpha^2 m_0 \int_{\S^1} |\partial_x^{k+1}(v^\alpha-w^\alpha)|^2\, dx
    -
    \alpha^2\int_{\S^1} \partial_x^{k+1}(v^\alpha-w^\alpha)\,  \partial_x^{k}\big((m(w^\alpha)-m_0) (v^\alpha-w^\alpha)_x\big)\, dx.
\end{split}
\end{equation*}
Applying the Cauchy--Schwarz inequality to the second term on the right-hand side, using that $H^k(\S^1)$ is a Banach algebra and Poincaré's inequality yield
\begin{equation*}
\begin{split}
    &-\alpha^2\int_{\S^1} \partial_x^{k+1}(v^\alpha-w^\alpha)\,  \partial_x^{k}\big((m(w^\alpha)-m_0) (v^\alpha-w^\alpha)_x\big)\, dx\\
    & \leq
    \alpha^2
    \|\partial_x^{k+1}(v^\alpha-w^\alpha)\|_{L_2}
    \|\partial_x^{k}\big((m(w^\alpha)-m_0) (v^\alpha-w^\alpha)_x\big)\|_{L_2}
    \\
    &\leq
    \alpha^2
    \|\partial_x^{k+1}(v^\alpha-w^\alpha)\|_{L_2}
    \|(m(w^\alpha)-m_0) (v^\alpha-w^\alpha)_x\|_{H^k}
    \\
    &\leq
    C\alpha^2
    \|\partial_x^{k+1}(v^\alpha-w^\alpha)\|_{L_2}
    \|m(w^\alpha)-m_0\|_{H^k} \|(v^\alpha-w^\alpha)_x\|_{H^k}
    \\
    &\leq
    C\alpha^2
    \|m(w^\alpha)-m_0\|_{H^k} \|\partial_x^{k+1}(v^\alpha-w^\alpha)\|_{L_2}^2.
\end{split}
\end{equation*}
Thus, for all $0 < t < T_\eps$ we have
\begin{equation} \label{eq:time_derivative_higher_energy}
    \frac{d}{dt} E_\alpha^k[v^\alpha,w^\alpha](t)
    \leq
    -\alpha^2 m_0 \int_{\S^1} |\partial_x^{k+1}(v^\alpha-w^\alpha)|^2\, dx
    +
    C\alpha^2
    \|m(w^\alpha)-m_0\|_{H^k} \|\partial_x^{k+1}(v^\alpha-w^\alpha)\|_{L_2}^2.
\end{equation}
We now show that $\|m(w^\alpha)-m_0\|_{H^k}$ is small in order to absorb the second term on the right-hand side into the first one.
Indeed, we have
\begin{equation}\label{eq:v-close-implies-w-close}
    \norm{w^\alpha - \hat{v}_0}_{H^k} \leq \norm{S_\alpha(v^\alpha-\hat{v}_0)}_{H^k} + \norm{S_\alpha \hat{v}_0 - \hat{v}_0}_{H^k}.
\end{equation}
By assumption, we have $\norm{v^\alpha - \hat{v}_0}_{H^k}\leq \eps$, hence it follows from Lemma \ref{lem:solution-operator-2nd} (ii) that $\norm{S_\alpha(v^\alpha-\hat{v}_0)}_{H^k}\leq C\eps$. On the other hand, we have $\norm{S_\alpha \hat{v}_0 - \hat{v}_0}_{H^k} \leq \tfrac{C}{\alpha-1}\hat{v}_0 \leq C\eps$ if $\alpha$ is large enough. Thus, $w^\alpha$ is $H^k$-close to the constant $\hat{v}_0$, too. Thence, multiple applications of the chain rule and the mean value theorem yield that also $m(w^\alpha)$ must be $H^k$-close to $m_0 = m(\hat{v}_0)$:
\begin{equation}\label{eq:m-almost-constant}
    \|m(w^\alpha)-m_0\|_{H^k}
    \leq
    C_{k,m_0} \eps
    \leq
    C_{k,m_0} \eps_0,
    \quad 0 < t < T_\eps.
\end{equation}
Inserting this into \eqref{eq:time_derivative_higher_energy} and choosing $\eps_0 > 0$ small enough, we end up with
\begin{equation*}
    \frac{d}{dt} E_\alpha^k[v^\alpha,w^\alpha](t)
    \leq
    -\frac{\alpha^2 m_0}{2} \int_{\S^1} |\partial_x^{k+1}(v^\alpha-w^\alpha)|^2\, dx,
    \quad 0 < t < T_\eps.
\end{equation*}
It remains to estimate the dissipation term on the right-hand side against the energy. This will allow us to apply Gronwall's inequality and thus to finally obtain exponential decay. 
Indeed, splitting $\alpha^2 = \frac{\alpha^2}{\alpha +1} + \frac{\alpha^3}{\alpha+1}$ and using $\eqref{eq:relaxed_system_conv}_2$, we can write for all $0 < t < T_\eps$
\begin{equation}\label{eq:higher-order_energy_dissipation}
\begin{split}
    \frac{d}{dt} E_\alpha^k[v^\alpha,w^\alpha](t)
    &\leq
    -\frac{\alpha^2 m_0}{2} \int_{\S^1} |\partial_x^{k+1}(v^\alpha-w^\alpha)|^2\, dx
    \\
    &=
    -\frac{\alpha^2 m_0}{2(\alpha +1)} \int_{\S^1} |\partial_x^{k+1}(v^\alpha-w^\alpha)|^2\, dx
    -
    \frac{\alpha^3 m_0}{2(\alpha+1)} \int_{\S^1} |\partial_x^{k+1}(v^\alpha-w^\alpha)|^2\, dx
    \\
    &=
    -\frac{\alpha^2 m_0}{2(\alpha +1)} \int_{\S^1} |\partial_x^{k+1}(v^\alpha-w^\alpha)|^2\, dx
    -
    \frac{\alpha m_0}{2(\alpha +1)}
    \int_{\S^1} |\partial_x^{k+3} w^\alpha + \partial_x^{k+1}w^\alpha|^2\, dx.
\end{split}
\end{equation}
Using Poincaré's inequality, the first term on the right-hand side can be estimated from below by
\begin{equation}\label{eq:dissipation_1}
    \frac{\alpha^2 m_0}{2(\alpha +1)} \int_{\S^1} |\partial_x^{k+1}(v^\alpha-w^\alpha)|^2\, dx
    \geq 
    \frac{\alpha m_0}{\alpha +1} 
    \left(\frac{\alpha}{2}\int_{\S^1} |\partial_x^{k}(v^\alpha-w^\alpha)|^2\, dx\right),
\end{equation}
while for the second term, we use the Fourier representation of $w^\alpha$, to obtain
\begin{equation} \label{eq:dissipation_2}
    \frac{\alpha m_0}{2(\alpha +1)}
    \int_{\S^1} |\partial_x^{k+3} w^\alpha + \partial_x^{k+1}w^\alpha|^2\, dx
    \geq
    \frac{\alpha m_0}{\alpha+1} \left(\frac{1}{2} \int_{\S^1} |\partial_x^{k+1}w^\alpha|^2 - |\partial_x^k w^\alpha|^2\, dx\right).
\end{equation}
Plugging \eqref{eq:dissipation_1} and \eqref{eq:dissipation_2} into \eqref{eq:higher-order_energy_dissipation} leads to
\begin{equation*}
    \frac{d}{dt} E_\alpha^k[v^\alpha,w^\alpha](t)
    \leq
    -\frac{\alpha m_0}{\alpha+1} E_\alpha^k[v^\alpha,w^\alpha](t),
    \quad 0 < t < T_\eps.
\end{equation*}
By Gronwall's inequality, we end up with
\begin{equation} \label{eq:exponential_decay}
    E_\alpha^k[v^\alpha,w^\alpha](t)
    \leq 
    E_\alpha^k[v_0,w_0^\alpha]\, \text{exp}\big({-\tfrac{\alpha m_0}{\alpha+1} t}\big),
    \quad 0 < t < T_\eps.
\end{equation}
Here, the initial value $w_0^\alpha$ for $w^\alpha$ is obtained by solving $\eqref{eq:relaxed_system_conv}_2$ with $v = v_0$, that is, by setting $w_0^\alpha \coloneq S_\alpha v_0$. By assumption, $v_0$ is $\delta$-close to $\hat{v}_0$ and we have seen in \eqref{eq:v-close-implies-w-close} that the same is true for $w^\alpha_0$. This implies
\begin{equation}\label{eq:init_energy_smallness}
    E_\alpha^k[v_0,w_0^\alpha] \leq C \delta
\end{equation}
and hence exponential decay of the higher-order energy.

By definition of the higher-order energy $E_\alpha^k$ and \eqref{eq:exponential_decay} we control all but the first Fourier modes of $w^\alpha$ (cf. the Fourier representation of the higher-order energy in \eqref{eq:higher-order-energy-fourier-representation}). Thanks to the term $\frac{\alpha}{2}\|\partial_x^k(v^\alpha-w^\alpha)\|_{L_2}^2$ in the energy, we then also control all but the first Fourier modes of $v^\alpha$.
In order to estimate the missing first Fourier modes, we use the fundamental theorem of calculus and write for all $0 < t < T_\eps$
\begin{equation*}
    \hat v_{\pm 1}(t)\,e^{\pm \ii x} = \hat v_{\pm 1}(0)\,e^{\pm\ii x} + \int_0^t P_{\pm 1}\big[\big(m(w^\alpha)\alpha (v^\alpha-w^\alpha)_x\big)_x\big]\, ds,
\end{equation*}
where $P_{\pm 1}$ denotes the projection onto the first Fourier modes, and
\begin{equation*}
\begin{aligned}
    \normb{P_{\pm 1}\big[\big(m(w^\alpha)\alpha (v^\alpha-w^\alpha)_x\big)_x\big]}_{L_2} \leq &\, \normb{P_{\pm 1}\big[\big((m(w^\alpha)-m_0)\alpha (v^\alpha-w^\alpha)_x\big)_x\big]}_{L_2}\\
    & + m_0 \normb{P_{\pm 1}\big[\alpha (v^\alpha-w^\alpha)_{xx}\big]}_{L_2}\\
    \leq & \,\normb{(m(w^\alpha)-m_0)\alpha (v^\alpha-w^\alpha)_x}_{H^1} + m_0 \normb{P_{\pm 1}\big[-w_{xxxx}^\alpha-w_{xx}^\alpha\big]}_{L_2}\\
    \leq & \,C\norm{m(w^\alpha)-m_0}_{H^1} \norm{\alpha\partial_x^2(v^\alpha-w^\alpha)}_{L_2}.
\end{aligned}
\end{equation*}
Here, we have used that $H^1(\S^1)$ is a Banach algebra, Poincaré's inequality and the fact that $w_{xxxx} + w_{xx}$ does not contain first Fourier modes. Thence, using \eqref{eq:m-almost-constant} and $\eqref{eq:relaxed_system_conv}_2$ we estimate
\begin{equation*}
\begin{split}
    \abs{\hat v_{\pm 1}(t)}
    &\leq
    \abs{\hat v_{\pm 1}(0)}
    +
    \int_0^t \normb{P_{\pm 1}\big[\big(m(w^\alpha)\alpha (v^\alpha-w^\alpha)_x\big)_x\big]}_{L_2}ds
    \\
    &\leq
    \abs{\hat v_{\pm 1}(0)}
    +
    C_{k,m_0}\eps\int_0^t \norm{w_{xxxx}^\alpha + w_{xx}^\alpha}_{L_2} ds
    \\
    &\leq
    \abs{\hat v_{\pm 1}(0)}
    +
    C_{k,m_0}\eps\int_0^t \sqrt{E_\alpha^3[v^\alpha,w^\alpha](s)}\, ds.
\end{split}
\end{equation*}
Since $\norm{v_0 - \hat{v}_0}_{H^k} < \delta$, we have in particular $\abs{\hat v_{\pm 1}(0)} \leq C\delta$, and inserting \eqref{eq:exponential_decay} and \eqref{eq:init_energy_smallness} we find
\begin{equation*}
    \abs{\hat v_{\pm 1}(t)} \leq C_{k,m_0}(\delta + \eps_0\sqrt{\delta}).
\end{equation*}
On the other hand, we also control Fourier modes $\geq 2$ of $v$ by means of the energy:
\begin{equation*}
    \norm{P_{\geq 2} v^\alpha(t)}_{H^k} \leq \norm{P_{\geq 2} (v^\alpha-w^\alpha)(t)}_{H^k} + \norm{P_{\geq 2} w^\alpha(t)}_{H^k} \leq C_k \sqrt{E_\alpha^{k}[v^\alpha,w^\alpha](t)} \leq C\sqrt{\delta}.
\end{equation*}
Hence, choosing $\delta = \delta(\eps)>0$ small enough, we finally conclude
\begin{equation*}
    \norm{v^\alpha(t) - \hat{v}_0}_{H^k} \leq C_{k,m_0}\sqrt{\delta} < \eps,
\end{equation*}
thus contradicting $T_\eps < \infty$.
This concludes the proof of global existence and uniqueness of smooth solutions.

It remains to prove the convergence \eqref{eq:first-fourier-convergence}. To this end, we estimate for $0\leq t_1\leq t_2$
\begin{equation*}
\begin{aligned}
    \norm{\hat v_{\pm 1}(t_1) - \hat v_{\pm 1}(t_2)}_{L^2}
    &\leq \int_{t_1}^{t_2}\normb{\tfrac{d}{ds}\hat v_{\pm 1}(s)}_{L^2}\, ds\\
    &\leq C\int_{t_1}^{t_2} \norm{m(w^\alpha(s))-m_0}_{H^1}\; \norm{\alpha\partial_x^2(v^\alpha(s)-w^\alpha(s))}_{L_2} ds\\
    &\leq C_{k,m_0}\eps\sqrt\delta\int_{t_1}^{t_2}\mathrm{exp}\big(-\tfrac{\alpha m_0}{\alpha+1}s\big) ds\\
    &\leq C_{k,m_0} \eps\sqrt\delta\, \mathrm{exp}\big(-\tfrac{\alpha m_0}{\alpha+1}t_1\big),
\end{aligned}
\end{equation*}
which converges to zero as $t_1\to\infty$. We point out that, in fact, the choice of norm on the first Fourier modes does not matter as the space of first Fourier modes is finite-dimensional. Thus, there exist $b_\alpha, c_\alpha\in \R$, such that
\begin{equation*}
    \lim_{t\to\infty} \big(\hat v_{1}(t)\ee^{\ii x} + \hat v_{-1}(t)\ee^{-\ii x}\big) = b_\alpha\cos(x) + c_\alpha\sin(x)
\end{equation*}
with the exponential estimate
\begin{equation*}
    \norm{P_{\pm 1}[v^\alpha(t)] - b_\alpha\cos(x) - c_\alpha\sin(x)} \leq C\eps\,\mathrm{exp}\big(-\tfrac{\alpha m_0}{\alpha+1}t\big).
\end{equation*}
This estimate is complemented by
\begin{equation}\label{eq:stability_decay}
    \norm{P_{\geq 2}[v^\alpha(t)]}_{H^k} \leq C\sqrt{E_\alpha^{k}[v^\alpha,w^\alpha](t)} \leq C\eps\, \mathrm{exp}\big({-\tfrac{\alpha m_0}{2(\alpha+1)} t}\big).
\end{equation}
Using the mass conservation, this concludes the proof of the exponential convergence for $v^\alpha$.
The estimate for $w^\alpha$ follows by applying the solution operator $S_\alpha$ to $v^\alpha$ and using its boundedness from Lemma \ref{lem:solution-operator-2nd} (ii).
\end{proof}

\begin{remark}\label{rem:decay-rate-uniform-in-alpha}
    Observe that, as $\alpha \to \infty$, the exponential rate of convergence obtained in
    \begin{equation*}
    E_\alpha^k[v^\alpha,w^\alpha](t)
    \leq 
    E_\alpha^k[v_0,w^\alpha_0]\, \mathrm{exp}\big({-\tfrac{\alpha m_0}{\alpha+1} t}\big),
    \quad t > 0,
\end{equation*}
approaches $m_0$.
\end{remark}

\medskip
\noindent\textbf{\textsc{Consistency of long-time behaviour in the limit $\alpha\to\infty$.}}
We have seen in the previous subsection that the long-time dynamics of both the relaxation system \eqref{eq:relaxed_system_conv} (see Theorem \ref{thm:stability-relaxed}) and the rimming-flow equation \eqref{eq:rimming-flow_conv} (see Theorem \ref{thm:stability-cited}) are very similar: Starting out from an initial state close enough to a positive constant steady state, solutions converge exponentially to a specific steady state that differs from the constant only by (small) first Fourier modes, or, equivalently, a linear combination of $\sin(x)$ and $\cos(x)$. Thus, we now fix a single initial value $u_0 = v_0 \in H^4(\S^1)$ for both for the rimming-flow problem as well as the relaxed problem. In particular, their mean values must coincide, i.e. $\hat{u}_0 = \hat{v}_0$. We assume that the initial value is close enough to its mean such that the assumptions of both Theorem \ref{thm:stability-cited} and Theorem \ref{thm:stability-relaxed} are fulfilled. The question of long-time consistency then reduces to whether we can show that $\lim_{\alpha\to\infty} v_\ast^\alpha = u_\ast$. This is answered affirmatively in the following theorem:

\begin{theorem}[Long-time consistency]\label{thm:long-time-consistency}
    Let $u_0 = v_0\in H^4(\S^1)$ be a positive initial value, i.e. $u_0(x) > 0$ for all $x\in\S^1$. Suppose that it satisfies the assumptions of both Theorem \ref{thm:stability-cited} and Theorem \ref{thm:stability-relaxed} (e.g. for $k=3$). Then, we have
    \begin{equation}\label{eq:long-time-consistency}
        \lim_{\alpha\to\infty} v_\ast^\alpha = \lim_{\alpha\to\infty}w_\ast^\alpha = u_\ast.
    \end{equation}
\end{theorem}

\begin{remark}
\begin{enumerate}[label=(\roman*)]
    \item Since $v_\ast^\alpha,w_\ast^\alpha,u_\ast \in \hat{u}_0 + P_{\pm 1} L^2(\S^1)$ are all elements of a two-dimensional affine space, it does not matter in which function space the convergence \eqref{eq:long-time-consistency} takes place. In fact, the convergence in any norm is simply equivalent to the convergence of the coefficients: $\lim_{\alpha\to\infty} b_\alpha = b$ and $\lim_{\alpha\to\infty} c_\alpha = c$.
    \item We briefly comment on the significance of the long-time consistency result in Theorem \ref{thm:long-time-consistency}: Depending on a given initial value, the rimming-flow equation chooses a specific steady state to which it converges in the long-time limit (cf. Theorem \ref{thm:stability-cited}). As it stands, the mechanism behind this choice is not clear. Thanks to the consistency result it will be possible to approach this question by further investigating the relaxation system (e.g. numerically) and possibly gain a better understanding of the dynamics in the original equation. However, this is outside the scope of this paper.
\end{enumerate}
\end{remark}

\begin{proof}
    It suffices to show $\lim_{\alpha\to\infty} v^\alpha_\ast = u_\ast$, the convergence of $w^\alpha_\ast$ to $u_\ast$ is then a consequence of $w^\alpha_\ast = S_\alpha v^\alpha_\ast$. Using that $u_\ast\in \hat u_0 + P_{\pm 1} L^2(\S^1)$, the fundamental theorem of calculus yields
    \begin{equation*}
        u_\ast = \hat u_0 + P_{\pm 1} u_\ast = \hat u_0 + P_{\pm 1}u_0 - \int_0^\infty P_{\pm 1}\left[m(u) (u_{xxx} + u_x)\right]_x\,dt.
    \end{equation*}
    Applying the same procedure to $v^\alpha$ yields
    \begin{equation*}
        v^\alpha_\ast = \hat v_0 + P_{\pm 1} v^\alpha_\ast = \hat v_0 + P_{\pm 1}v_0 - \int_0^\infty P_{\pm 1}\left[m(w^\alpha) (w^\alpha_{xxx} + w^{\alpha}_x)\right]_x\,dt.
    \end{equation*}
    Since $u_0 = v_0$ and $\hat u_0 = \hat v_0$, we need to show
    \begin{equation*}
        \lim_{\alpha\to\infty} \int_0^\infty P_{\pm 1}\left[m(w^\alpha) (w^\alpha_{xxx} + w^{\alpha}_x)\right]_x\,dt = \int_0^\infty P_{\pm 1}\left[m(u) (u_{xxx} + u_x)\right]_x\,dt.
    \end{equation*}
    To this end, let $\tilde\eps>0$ be arbitrary. We use $\tilde\eps$ instead of $\eps$ to prevent confusion with the $\eps$ from Theorem \ref{thm:stability-relaxed}. We split the first integral at a time $\tau_{\tilde\eps}$ to be chosen later depending on $\tilde\eps$:
    \begin{equation*}
    \begin{aligned}
        \int_0^\infty P_{\pm 1}\left[m(w^\alpha) (w^\alpha_{xxx} + w^{\alpha}_x)\right]_x\,dt
        &
        = \int_0^{\tau_{\tilde\eps}} P_{\pm 1}\left[m(w^\alpha) (w^\alpha_{xxx} + w^{\alpha}_x)\right]_x\,dt\\
        &
       + \int_{\tau_{\tilde\eps}}^\infty P_{\pm 1}\left[m(w^\alpha) (w^\alpha_{xxx} + w^{\alpha}_x)\right]_x\,dt.
    \end{aligned}
    \end{equation*}
    The second integral on the right-hand side is arbitrarily small, provided $\tau_\eps$ is large enough. Indeed, we can estimate
    \begin{equation}
    \begin{aligned}
        \normBB{\int_{\tau_{\tilde\eps}}^\infty P_{\pm 1}\left[m(w^\alpha) (w^\alpha_{xxx} + w^{\alpha}_x)\right]_x\,dt}_{L^2}
        &\leq \int_{\tau_{\tilde\eps}}^\infty \normb{P_{\pm 1}\left[m(w^\alpha) (w^\alpha_{xxx} + w^{\alpha}_x)\right]_x} \,dt\\
        &\leq C \int_{\tau_{\tilde\eps}}^\infty \norm{m(w^\alpha)}_{H^1}\,\norm{w^\alpha_{xxx} + w^\alpha_x}_{H^1}\,dt\\
        &\leq C\int_{\tau_{\tilde\eps}}^\infty\norm{P_{\geq 2} w^\alpha}_{H^4}\\
        &= C\int_{\tau_{\tilde\eps}}^\infty\norm{P_{\geq 2} \,S_\alpha v^\alpha}_{H^4}\,dt.
    \end{aligned}
    \end{equation}
    Hence, by Lemma \ref{lem:solution-operator-2nd} (ii) and the exponential decay estimate \eqref{eq:stability_decay}, we find
    \begin{equation*}
        \int_{\tau_{\tilde\eps}}^\infty\norm{P_{\geq 2} \,S_\alpha v^\alpha}_{H^4}\,dt \leq C_{m_0} \,\eps \,\exp\big(-\tfrac{\alpha m_0}{2(\alpha + 1)}\tau_{\tilde\eps}\big) < \tilde\eps,
    \end{equation*}
    provided $\tau_{\tilde\eps}$ is large enough. Note that $\tau_{\tilde\eps}$ may be chosen independently of $\alpha$, thanks to Remark \ref{rem:decay-rate-uniform-in-alpha}. Therefore, it remains to prove that
    \begin{equation}\label{eq:consistency-initial-integral-convergence}
        \lim_{\alpha\to\infty} \int_0^{\tau_{\tilde\eps}} P_{\pm 1}\left[m(w^\alpha) (w^\alpha_{xxx} + w^{\alpha}_x)\right]_x\,dt = \int_0^{\tau_{\tilde\eps}} P_{\pm 1}\left[m(u) (u_{xxx} + u_x)\right]_x\,dt.
    \end{equation}
    From the exponential decay \eqref{eq:exponential_decay} with $k\geq 3$ and the non-negativity of $w^\alpha$, it follows in particular that $w^\alpha$ is uniformly bounded in $H^4(\S^1)$:
    \begin{equation*}
        \norm{w^\alpha(t)}_{H^4} \leq C, \quad 0 \leq t < \infty.
    \end{equation*}
    Moreover, due to Lemma \ref{lem:hoelder-continuity}
    $w^\alpha$ is bounded in $C^{1/8}\big([0,\tau_{\tilde\eps}];L^2(\S^1)\big)$. Hence, an interpolation argument shows that $w^\alpha$ must, in fact, be bounded in $C^{1/32}\big([0,\tau_{\tilde\eps}];H^3(\S^1)\big)$:
    \begin{equation*}
        \norm{w^\alpha(t_2) - w^\alpha(t_1)}_{H^3} \leq \norm{w^\alpha(t_2) - w^\alpha(t_1)}_{L^2}^{1/4} \ \norm{w^\alpha(t_2) - w^\alpha(t_1)}_{H^4}^{3/4} \leq C\abs{t_2 - t_1}^{1/32}.
    \end{equation*}
    An Aubin--Lions--Simon result, cf. \cite[Theorem 5]{simon1986compact}, implies compactness of the embedding
    \begin{equation*}
        C^{1/32}\big([0,\tau_{\tilde\eps}];H^3(\S^1)\big)\cap L^\infty\big((0,\tau_{\tilde\eps});H^4(\S^1)\big) \overset{c}{\hookrightarrow} C\big([0,\tau_{\tilde\eps}];H^3(\S^1)\big)
    \end{equation*}
    As a consequence, every sequence $(w^{\alpha_n})_{n\in\N}$ must have a subsequence $(w^{\alpha_{n_k}})_{k\in\N}$ that converges strongly in $C\big([0,\tau_{\tilde\eps}];H^3(\S^1)\big)$. The limit of every such subsequence is a \emph{positive} weak solution of the rimming-flow problem by Theorem \ref{thm:convergence} and thus unique. Therefore, we conclude that
    \begin{equation*}
        \lim_{\alpha\to\infty} w^\alpha = u \quad\text{in }C\big([0,\tau_{\tilde\eps}];H^3(\S^1)\big),
    \end{equation*}
    whence it follows in particular that
    \begin{equation*}
        \lim_{\alpha\to\infty} P_{\pm 1}\left[m(w^\alpha) (w^\alpha_{xxx} + w^{\alpha}_x)\right] = P_{\pm 1}\left[m(u) (u_{xxx} + u_x)\right]
    \end{equation*}
    in $C\big([0,\tau_{\tilde\eps}]; P_{\pm 1} L^2(\S^1)\big)$. Now, for an arbitrary function $f\in P_{\pm 1} L^2(\S^1)$, we have also $f_x\in P_{\pm 1} L^2(\S^1)$ with identical norm $\norm{f}_{L^2} = \norm{f_x}_{L^2}$. This allows us to conclude
    \begin{equation*}
        \lim_{\alpha\to\infty} P_{\pm 1}\left[m(w^\alpha) (w^\alpha_{xxx} + w^{\alpha}_x)\right]_x = P_{\pm 1}\left[m(u) (u_{xxx} + u_x)\right]_x
    \end{equation*}
    uniformly in $[0,\tau_{\tilde\eps}]$, and \eqref{eq:consistency-initial-integral-convergence} follows.
\end{proof}


\section*{Statements and Declarations}

\noindent\textsc{Conflict of Interest. }
The authors have no competing interests to declare that are relevant to the content of this article.
\medskip

\noindent\textsc{Data Availability. }
Data sharing is not applicable to this article as no datasets were generated or analysed during the current study.
\medskip

\noindent\textsc{Funding. }
The authors were funded by Deutsche Forschungsgemeinschaft (DFG,
German Research Foundation) under Germany´s
Excellence Strategy – EXC 2075 – 390740016.
J. Joussen und C. Lienstromberg acknowledge support by the Deutsche Forschungsgemeinschaft (DFG, German Research Foundation) -- project number 545145736 -- through the scientific network \emph{(In)stability Phenomena in Asymptotic Models in Fluid Dynamics}.
\medskip

\noindent\textsc{Use of Generative Artificial Intelligence. }
No generative artificial intelligence tools were used, neither for the mathematical results nor for the writing of the manuscript.
\medskip

\printbibliography

@article{jlv26,
  title={Two-time scale dynamics of solutions to a rimming-flow equation},
  author={Joussen, Juri and Lienstromberg, Christina and Vel{\'a}zquez, Juan JL},
  journal={Pure Appl. Anal.},
  volume={8},
  number={3},
  pages={681--719},
  year={2026},
  publisher={Mathematical Sciences Publishers}
}

@article{pernas_castano_analysis_2020,
	title = {Analysis of a thin film approximation for two-fluid Taylor-Couette flows},
	volume = {269},
	issn = {0022-0396},
	doi = {10.1016/j.jde.2019.12.005},
	pages = {377--419},
	number = {1},
	fjournal = {Journal of Differential Equations},
	journal = {J. Differ. Equ.},
	author = {Pernas Castaño, Tania and Velázquez, Juan J. L.},
	urldate = {2024-02-01},
	date = {2020-06-15},
}

@article{karabash_multi-parameter_2024,
  title={Multi-parameter Hopf bifurcations of rimming flows},
  author={Karabash, Illya M and Lienstromberg, Christina and Vel{\'a}zquez, Juan JL},
  journal={J. Differ. Equ.},
  volume={430},
  pages={113182},
  year={2025},
  publisher={Elsevier}
}

@article{lienstromberg_analysis_2022,
	title = {Analysis of a Two-Fluid Taylor–Couette Flow with One Non-Newtonian Fluid},
	volume = {32},
	issn = {1432-1467},
	doi = {10.1007/s00332-021-09750-0},
	pages = {24},
	number = {2},
	fjournal = {Journal of Nonlinear Science},
	journal = {J. Nonlinear Sci.},
	author = {Lienstromberg, Christina and Pernas-Castaño, Tania and Velázquez, Juan J. L.},
	urldate = {2023-01-11},
	date = {2022-02-12},
	langid = {english},
}

@article{Seiden_Thomas_2011,
  title = {Complexity, segregation, and pattern formation in rotating-drum flows},
  author = {Seiden, G. and Thomas, P. J.},
  journal = {Rev. Mod. Phys.},
  volume = {83},
  issue = {4},
  pages = {1323--1365},
  numpages = {0},
  year = {2011},
  month = {11},
  publisher = {American Physical Society},
  doi = {10.1103/RevModPhys.83.1323},
}

@article {Moffatt_1977,
    AUTHOR = {Moffatt, H. K.},
     TITLE = {Behaviour of a viscous film on the outer surface of a rotating cylinder},
   JOURNAL = {J. de Mécanique},
  FJOURNAL = {Journal de Mécanique},
    VOLUME = {16},
      YEAR = {1977},
    NUMBER = {5},
     PAGES = {651--673},
      ISSN = {},
   MRCLASS = {},
  MRNUMBER = {},
       DOI = {},
}

@article{BF1990,
title = {Higher order nonlinear degenerate parabolic equations},
journal = {J. Differ. Equ.},
volume = {83},
number = {1},
pages = {179-206},
year = {1990},
issn = {0022-0396},
doi = {https://doi.org/10.1016/0022-0396(90)90074-Y},
author = {Francisco Bernis and Avner Friedman}
}

@article {LV2024,
    AUTHOR = {Lienstromberg, Christina and Vel\'azquez, Juan J. L.},
     TITLE = {Long-time asymptotics and regularity estimates for weak
              solutions to a doubly degenerate thin-film equation in the
              {T}aylor-{C}ouette setting},
   JOURNAL = {Pure Appl. Anal.},
  FJOURNAL = {Pure and Applied Analysis},
    VOLUME = {6},
      YEAR = {2024},
    NUMBER = {1},
     PAGES = {187--236},
      ISSN = {2578-5885,2578-5893},
   MRCLASS = {35B40 (35K35 35K65 35Q35 76A05 76A20)},
  MRNUMBER = {4707707},
MRREVIEWER = {Hujun\ Yang},
       DOI = {10.2140/paa.2024.6.187},
}

@article {BBdP1995,
    AUTHOR = {Beretta, Elena and Bertsch, Michiel and Dal Passo, Roberta},
     TITLE = {Nonnegative solutions of a fourth-order nonlinear degenerate parabolic equation},
   JOURNAL = {Arch. Rational Mech. Anal.},
  FJOURNAL = {Archive for Rational Mechanics and Analysis},
    VOLUME = {129},
      YEAR = {1995},
    NUMBER = {2},
     PAGES = {175--200},
      ISSN = {0003-9527},
   MRCLASS = {35K65 (35B05 35K55)},
  MRNUMBER = {1328475},
MRREVIEWER = {Minkyu\ Kwak},
       DOI = {10.1007/BF00379920},
}

@article{Throne_Gianchandani_1980,
author = {Throne, J. L. and Gianchandani, J.},
title = {Reactive rotational molding},
journal = {Polym. Eng. Sci.},
volume = {20},
number = {13},
pages = {899-919},
doi = {https://doi.org/10.1002/pen.760201309},
year = {1980}
}

@article{white1956residual,
  title={Residual condensate, condensate behavior, and siphoning in paper driers},
  author={White, ROBERT E},
  journal={Tech. Assoc. Pulp Paper Ind},
  volume={39},
  pages={228},
  year={1956}
}

@article{yih1960instability,
  title={Instability of a rotating liquid film with a free surface},
  author={Yih, Chia-Shun},
  journal={Proc. R. Soc. Lond. A. Math. Phys. Sci.},
  volume={258},
  number={1292},
  pages={63--89},
  year={1960},
  publisher={The Royal Society London}
}

@book{lunardi_1995,
 author = {Lunardi, Alessandra},
 title = {Analytic semigroups and optimal regularity in parabolic problems},
 fseries = {Progress in Nonlinear Differential Equations and Their Applications},
 series = {Prog. Nonlinear Differ. Equ. Appl.},
 issn = {1421-1750},
 volume = {16},
 isbn = {3-7643-5172-1},
 year = {1995},
 publisher = {Basel: Birkh{\"a}user},
 language = {English},
 zbMATH = {732330},
 Zbl = {0816.35001}
}

@article{pukhnachov2005capillary,
  title={Capillary/gravity film flows on the surface of a rotating cylinder},
  author={Pukhnachov, VV},
  fjournal={Journal of Mathematical Sciences},
  journal={J. Math. Sci.},
  volume={130},
  number={4},
  pages={4871--4883},
  year={2005},
  publisher={Springer}
}

@article{benilov2005does,
  title={Does surface tension stabilize liquid films inside a rotating horizontal cylinder?},
  author={Benilov, ES and Kopteva, N and O'Brien, SBG},
  fjournal={Quarterly Journal of Mechanics and Applied Mathematics},
  journal={Q. J. Mech. Appl. Math.},
  volume={58},
  number={2},
  pages={185--200},
  year={2005},
  publisher={OUP}
}

@article{burchard2012convergence,
  title={Convergence to equilibrium for a thin-film equation on a cylindrical surface},
  author={Burchard, Almut and Chugunova, Marina and Stephens, Benjamin K},
  journal={Commun. Partial Differ. Equ.},
  volume={37},
  number={4},
  pages={585--609},
  year={2012},
  publisher={Taylor \& Francis}
}

@article{chugunova2010nonnegative,
  title={Nonnegative solutions for a long-wave unstable thin film equation with convection},
  author={Chugunova, Marina and Pugh, Mary C and Taranets, Roman M},
  JOURNAL = {SIAM J. Math. Anal.},
  FJOURNAL={SIAM Journal on Mathematical Analysis},
  volume={42},
  number={4},
  pages={1826--1853},
  year={2010},
  publisher={SIAM}
}

@article{jllv2026,
      title={Analysis of a three-dimensional fluid flow in rotating cylinders}, 
      author={Juri Joussen and Janne Laudien and Christina Lienstromberg and Juan J. L. Velázquez},
      year={2026},
      eprint={2605.10305},
      archivePrefix={arXiv},
      primaryClass={math.AP},
      url={https://arxiv.org/abs/2605.10305}, 
}

@article{chugunova2009nature,
  title={On the nature of ill-posedness of the forward-backward heat equation},
  author={Chugunova, Marina and Karabash, Illya M and Pyatkov, Sergei G},
  journal={Integral Equ. Oper. Theory},
  volume={65},
  number={3},
  pages={319--344},
  year={2009},
  publisher={Springer}
}

@article{lopes2018multiple,
  title={On the multiple solutions of coating and rimming flows on rotating cylinders},
  author={Lopes, Andr{\'e} v B and Thiele, Uwe and Hazel, Andrew L},
  journal={J. Fluid Mech.},
  volume={835},
  pages={540--574},
  year={2018},
  publisher={Cambridge University Press}
}

@article{benilov2008steady,
  title={Steady rimming flows with surface tension},
  author={Benilov, ES and Benilov, MS and Kopteva, Natalia},
  journal={J. Fluid Mech.},
  volume={597},
  pages={91--118},
  year={2008},
  publisher={Cambridge University Press}
}

@article{pukhnachev1978motion,
  title={Motion of a liquid film on the surface of a rotating cylinder in a gravitational field},
  author={Pukhnachev, VV},
  journal={PMTF Zhurnal Prikladnoi Mekhaniki i Tekhnicheskoi Fiziki},
  volume={18},
  number={3},
  pages={78--88},
  year={1978}
}

@article{ashmore_hosoi_stone_2003,
 author = {Ashmore, J. and Hosoi, A. E. and Stone, H. A.},
 title = {The effect of surface tension on rimming flows in a partially filled rotating cylinder},
 fjournal = {Journal of Fluid Mechanics},
 journal = {J. Fluid Mech.},
 issn = {0022-1120},
 volume = {479},
 pages = {65--98},
 year = {2003},
 language = {English},
 doi = {10.1017/S0022112002003312},
 zbMATH = {2065441},
 Zbl = {1163.76350}
}

@article{thiele2011depinning,
  title={On the depinning of a drop of partially wetting liquid on a rotating cylinder},
  author={Thiele, Uwe},
  journal={J. Fluid Mech.},
  volume={671},
  pages={121--136},
  year={2011},
  publisher={Cambridge University Press}
}

@book{evans_partial_2010,
	address = {Providence, RI},
	edition = {2nd ed.},
	series = {Grad. {Stud}. {Math}.},
	title = {Partial differential equations},
	volume = {19},
	isbn = {978-0-8218-4974-3},
	issn = {1065-7338},
	language = {English},
	publisher = {American Mathematical Society (AMS)},
	author = {Evans, Lawrence C.},
	year = {2010},
}

@article{cornalba2025h,
  title={Hölder continuity of weak solutions to the thin-film equation in $ d= 2$},
  author={Cornalba, Federico and Fischer, Julian and Kokavcov{\'a}, Erika Maringov{\'a}},
  journal={arXiv preprint arXiv:2512.24809},
  year={2025}
}

@book{amann_linear_1995,
	address = {Basel},
	title = {Linear and {Quasilinear} {Parabolic} {Problems}: {Volume} {I}: {Abstract} {Linear} {Theory}},
	copyright = {http://www.springer.com/tdm},
	isbn = {978-3-0348-9950-5 978-3-0348-9221-6},
	url = {http://link.springer.com/10.1007/978-3-0348-9221-6},
	doi = {10.1007/978-3-0348-9221-6},
	language = {en},
	urldate = {2026-07-30},
	publisher = {Birkhäuser},
	author = {Amann, Herbert},
	year = {1995},
}

@article{cuesta_self-similar_2018,
	title = {Self-{Similar} {Lifting} and {Persistent} {Touch}-{Down} {Points} in the {Thin}-{Film} {Equation}},
	volume = {50},
	issn = {0036-1410, 1095-7154},
	url = {https://epubs.siam.org/doi/10.1137/17M1141485},
	doi = {10.1137/17M1141485},
	language = {en},
	number = {2},
	urldate = {2026-10-05},
	journal = {SIAM J. Math. Anal.},
	author = {Cuesta, Carlota M. and Knüpfer, Hans and Velázquez, Juan J. L.},
	month = jan,
	year = {2018},
	pages = {1900--1917},
}

@article{king_moving_2001,
	title = {Moving boundary problems and non-uniqueness for the thin film equation},
	volume = {12},
	url = {https://www.cambridge.org/core/journals/european-journal-of-applied-mathematics/article/moving-boundary-problems-and-nonuniqueness-for-the-thin-film-equation/E1EA284D6D9080E9856D990EDEEBAB1E},
	number = {3},
	urldate = {2026-10-05},
	journal = {Eur. J. Appl. Math.},
	publisher = {Cambridge University Press},
	author = {King, J. R. and Bowen, M.},
	year = {2001},
	pages = {321--356},
}

@article{simon1986compact,
  title={Compact sets in the space $L^p (O, T; B)$},
  author={Simon, Jacques},
  journal={Ann. Mat. Pura Appl.},
  volume={146},
  number={1},
  pages={65--96},
  year={1986},
  publisher={Springer}
}

@article{orr1978rimming,
  title={Rimming flow: numerical simulation of steady, viscous, free-surface flow with surface tension},
  author={Orr, FM and Scriven, LE},
  journal={J. Fluid Mech.},
  volume={84},
  number={1},
  pages={145--165},
  year={1978},
  publisher={Cambridge University Press}
}

@article{corli2012singular,
  title={Singular limits for a parabolic--elliptic regularization of scalar conservation laws},
  author={Corli, Andrea and Rohde, Christian},
  journal={J. Differ. Equ.},
  volume={253},
  number={5},
  pages={1399--1421},
  year={2012},
  publisher={Elsevier}
}

@article{corli2014parabolic,
  title={Parabolic approximations of diffusive--dispersive equations},
  author={Corli, Andrea and Rohde, Christian and Schleper, Veronika},
  journal={J. Math. Anal. Appl.},
  volume={414},
  number={2},
  pages={773--798},
  year={2014},
  publisher={Elsevier}
}

@article{barthwal2026hyperbolic,
  title={On hyperbolic approximations for a class of dispersive and diffusive-dispersive equations},
  author={Barthwal, Rahul and Dhaouadi, Firas and Rohde, Christian},
  journal={J. Nonlinear Sci.},
  volume={36},
  number={5},
  pages={106},
  year={2026},
  publisher={Springer}
}

@article{barthwal2026energy,
  title={Energy consistent hyperbolic approximation for a class of fourth-order partial differential equations},
  author={Barthwal, Rahul and Dhaouadi, Firas and Rohde, Christian},
  journal={arXiv preprint arXiv:2608.09478},
  year={2026}
}

@article{dhaouadi2025first,
    AUTHOR = {Dhaouadi, Firas and Dumbser, Michael and Gavrilyuk, Sergey},
     TITLE = {A first-order hyperbolic reformulation of the
              {C}ahn-{H}illiard equation},
   JOURNAL = {Proc. A},
  FJOURNAL = {Proceedings A},
    VOLUME = {481},
      YEAR = {2025},
    NUMBER = {2312},
     PAGES = {Paper No. 20240606, 24},
      ISSN = {1364-5021,1471-2946},
   MRCLASS = {35L40 (65M06)},
  MRNUMBER = {4926248},
MRREVIEWER = {Luigi\ Rarit\`a},
}

@article{taylor1923stability,
  title={Stability of a viscous liquid contained between two rotating cylinders},
  author={Taylor, Geoffrey Ingram},
  journal={Philos. Trans. R. Soc. London. Ser. A},
  volume={223},
  pages={289--343},
  year={1923},
  publisher={JSTOR}
}

@article{mallock1896experiments,
  title={Experiments on fluid viscosity},
  author={Mallock, Arnulph},
  journal={Philosophical Transactions of the Royal Society of London. Series A, Containing Papers of a Mathematical or Physical Character},
  volume={187},
  pages={41--56},
  year={1896},
  publisher={JSTOR}
}

@article{couette1890distinction,
  title={Distinction de deux r{\'e}gimes dans le mouvement des fluides},
  author={Couette, M},
  journal={J. Phys. Theor. Appl.},
  volume={9},
  number={1},
  pages={414--424},
  year={1890}
}

@article{ansini2004doubly,
  title={Doubly nonlinear thin-film equations in one space dimension},
  author={Ansini, Lidia and Giacomelli, Lorenzo},
  JOURNAL = {Arch. Rational Mech. Anal.},
  volume={173},
  number={1},
  pages={89--131},
  year={2004},
  publisher={Springer}
}

@article{perthame_relaxation_2021,
	title = {Relaxation of the {Cahn}–{Hilliard} equation with singular single-well potential and degenerate mobility},
	volume = {32},
	issn = {0956-7925, 1469-4425},
	doi = {10.1017/S0956792520000054},
	language = {en},
	number = {1},
	journal = {European J. Appl. Math.},
	author = {Perthame, Benoît and Poulain, Alexandre},
	month = feb,
	year = {2021},
	pages = {89--112},
}

\end{document}